\documentclass[12pt,a4paper]{amsart}
\usepackage[utf8]{inputenc}
\usepackage[english]{babel}
\usepackage[T1]{fontenc}
\usepackage{amsmath}
\usepackage{amsthm}
\usepackage{amsfonts}
\usepackage{amssymb}
\usepackage{multicol}
\usepackage{graphicx}

\usepackage[hidelinks]{hyperref}

\usepackage[all,cmtip]{xy}

\usepackage[left=3cm,right=3cm,top=3cm,bottom=3cm]{geometry}
\usepackage{pgf,tikz-cd}
\usetikzlibrary{decorations.pathreplacing}
\usetikzlibrary{arrows}
\usepackage{scalerel}[2016/12/29]

\usepackage{xcolor}
\definecolor{lime}{HTML}{A6CE39}
\DeclareRobustCommand{\orcidicon}{%
	\begin{tikzpicture}
	\draw[lime, fill=lime] (0,0) 
	circle [radius=0.16] 
	node[white] {{\fontfamily{qag}\selectfont \tiny ID}};
	\draw[white, fill=white] (-0.0625,0.095) 
	circle [radius=0.007];
	\end{tikzpicture}
	\hspace{-2mm}
}

\foreach \x in {A, ..., Z}{%
	\expandafter\xdef\csname orcid\x\endcsname{\noexpand\href{https://orcid.org/\csname orcidauthor\x\endcsname}{\noexpand\orcidicon}}
}

\DeclareMathOperator{\val}{val}
\DeclareMathOperator{\rv}{rv}
\DeclareMathOperator{\RV}{RV}
\DeclareMathOperator{\supp}{supp}

\title{On the hyperfields associated to valued fields}

\author[A.\ Linzi]{Alessandro Linzi\orcidA{}}
\address{Center for Information Technologies and Applied Mathematics\\ 
University of Nova Gorica\\ 
Slovenia}
\email{alessandro.linzi@ung.si}
\thanks{}

\author[P.\ Touchard]{Pierre Touchard\orcidB{}}
\address{Dipartimento di matematica e fisica, Università degli Studi della campania \lq\lq Luigi Vanvitelli\rq\rq.}
\email{pierre.pa.touchard@gmail.com}
\thanks{PT was supported by a grant of the University of Campania ‘Luigi Vanvitelli’ in the framework of V:ALERE 2019 (GoAL project)}

\theoremstyle{plain}
\newtheorem{theorem}{Theorem}[section]

\newtheorem{lemma}[theorem]{Lemma}
\newtheorem{proposition}[theorem]{Proposition}
\newtheorem{claim}{Claim}

\newtheorem{fact}[theorem]{Fact}

\theoremstyle{definition}
\newtheorem{definition}[theorem]{Definition}
\newtheorem{example}[theorem]{Example}
\newtheorem*{notation}{Notation}
\newtheorem{question}{Question}

\theoremstyle{remark}
\newtheorem{remark}[theorem]{Remark}

\begin{document}

\begin{abstract}
    One can associate to a valued field an inverse system of valued hyperfields $(\mathcal{H}_i)_{i \in I}$ in a natural way.  We investigate when, conversely, such a system arise from a valued field.  First, we extend a result of Krasner by showing that the inverse limit of certain systems are stringent valued hyperfields. Secondly, we describe a Hahn-like construction which yields a henselian valued field from a stringent valued hyperfield. In addition, we provide an axiomatisation of the theory of stringent valued hyperfields in a language consisting of two binary function symbols $\oplus$ and $\cdot$ and two constant symbols $\textbf{0}$ and $\textbf{1}$.
\end{abstract}

\maketitle

\section*{Introduction}

In the study of valued fields, the residue field and the value group play an important role. Since the pioneering work of  Ax, Kochen and Ershov (\cite{AK65A,AK65B,AK66,Ers65A,Ers65B,Ers65C,Ers65D}),  model theoretic properties of valued fields have been investigated and characterised at the level of residue fields and value groups. However, such reduction can be difficult to obtain, especially in the mixed characteristic case. Subsequent work has revealed then the advantage of considering intermediate structures, generalising both the value group and the residue field: let us mention the \textit{amc-structures} in the work of F.-V.\ Kuhlmann \cite{Kuh94}, and the \textit{mixed-structures} in the work of Basarab \cite{Bas91}. These ultimately lead Flenner in \cite{Fle11} to the definition of RV\emph{-structures}, which are currently used in the literature for various purposes  (see for instance \cite{CHR22,HM21,HK21}). The $\RV$-structures associated to a valued field $\mathcal{K}$ of mixed characteristic $(0,p)$ form a system $(\RV_n(\mathcal{K}))_{n\in \mathbb{N}}$ where, essentially, $\RV_n(\mathcal{K})$ is the multiplicative quotient group $K^\times/1+\mathfrak{m}_{n}$ equipped with some extra structure. Here, $\mathfrak{m}_n$ denotes the ideal of the valuation ring of $\mathcal{K}$ consisting of all those elements of value grater than the value of $p^n$. Flenner's main result  in \cite{Fle11}  states that a mixed characteristic henselian valued field $\mathcal{K}$ eliminates quantifiers relative to $(\RV_n(\mathcal{K}))_n$. He then reduced the problem of decidability of $\mathcal{K}$ to that of $(\RV_n(\mathcal{K}))_n$.

This result of quantifier elimination and the reduction method resulting from it, is used in a number of recent papers (e.g. \cite{ACGZ22,HM21,HK21}) and shows the intrinsic importance of the RV-structures: they capture more information than the value group and the residue field, and they are in many contexts simpler to analyse than the valued field to which they are associated. One can observe that in the above mentioned papers, the $\RV$-structures are not understood as algebraic structures on their own; but rather as an intermediate step before a reduction to the value group and the residue field.\par
Nevertheless, other authors consider $\RV$-structures independently of the valued field to which they are associated and as objects of study on their own. In this setting, $\RV$-structures are known as (an instance of) \textit{valued (Krasner) hyperfields}. A hyperfield $\mathcal{H}$ is a field where the (standard) additive operation is replaced by a \textit{multivalued operation}, that is, for all elements $a,b$ of $\mathcal{H}$ the sum $a+b$ is a subset of $\mathcal{H}$. A valued hyperfield is a hyperfield equipped with a valuation map (see Definition \ref{DefinitionValuationHyperfields}). More generally, these objects are instances of structures which are of interest for a well established research area known as hypercompositional algebra, initiated by Marty with the definition of the hypergroup in \cite{Mar34} (see also \cite{Gla,Mas21}).

The definition of valued hyperfield was first given in the work \cite{Kra57} of Krasner. The interested reader can consult the introduction of \cite{Vuk18} for a summary on earlier work of Krasner. For recent developments in the research on hyperfields see for instance \cite{KLS22,BS21}. We note that hyperfields have recently played a central role in an interesting paper of Junguk Lee on the model theory of valued fields (see \cite[Theorem 5.8]{Jun20}).

The main aim of our article is to give an answer to the following: 
\begin{question}\label{question}
Given a system of valued hyperfields $(\mathcal{H}_n)_n$, can we find a valued field $\mathcal{K}$ from which such a system arises, i.e., such that $\mathcal{H}_n=\RV_n(\mathcal{K})$ for all $n$?
\end{question}
After making this question more precise, we answer it positively and give a canonical construction of such a valued field. While dealing with the mentioned systems, we found it useful to combine and analyse  both points of view described above: that of hyperfields and that of $\RV$-structures.\par 
The paper is organized as follows. In Section \ref{SectionIsometricSequencesHyperfields}, after some preliminaries, we extend a result of Krasner by showing that the inverse limit of some isometric systems of valued hyperfields is a stringent (in the sense of  \cite{BS21}) valued hyperfield. In Section \ref{SectionAxiomatisationStringentHyperfields}, we give a useful axiomatisation of stringent valued hyperfields. Then, in Section \ref{SectionHahnProductStingentHyperfield}, we provide a Hahn-like construction which associates a (henselian) valued field to any stringent valued hyperfield. Finally, in Section \ref{SectionMainResult}, we state and prove our main result and briefly discuss how it answers Question \ref{question}.

\section*{Aknowledgements}
The authors would like to thank Franz-Viktor Kuhlmann for suggesting this topic and for all the fruitful discussions. Many thanks also to Paola D'Aquino for her advice and support.
\section{Isometric systems of valued hyperfields}\label{SectionIsometricSequencesHyperfields}

\subsection{Hyperfields and valuations}

Let $H$ be a nonempty set and denote by $\mathcal{P}^*(H)$ the family of all non-empty subsets of $H$. A \emph{hyperoperation} $+$ on $H$ is a function which associates with every pair $(x,y) \in H \times H$ an element of $\mathcal{P}^*(H)$, denoted by $x+y$. For $x \in H$ and $A,B\subseteq H$ we set 
\begin{equation}\label{A+B}
A+B=\bigcup_{a\in A,b\in B} a+b,
\end{equation}
$A + x = A + \lbrace x \rbrace$ and $x +A = \lbrace x \rbrace + A$. Note that if $A$ or $B$ is empty, then so is $A+B$.

The notion of hypergroup was defined by Marty in \cite{Mar34} to be a non-empty set $H$ with an associative hyperoperation $+$ (see Axiom (CH1) below) such that $x + H = H + x = H$ for all $x \in H$. Let us observe that in a hypergroup, sums are automatically non-empty. Indeed, suppose that $x+y=\emptyset$ for some $x,y\in H$. Then
\[
H=x+H=x+(y+H)=(x+y)+H=\emptyset+H=\emptyset,
\]
which is excluded.

The following special class of hypergroups will be of interest for us.

\begin{definition} \label{hypergp}
A \emph{canonical hypergroup} is a tuple $(H,+,0)$, where $(H, +)$ is a non-empty set with a hyperoperation and $0$ is an element of $H$ such that the following axioms hold:
\begin{itemize}
\item[(CH1)] (associativity) the hyperoperation $+$ is associative, i.e., $(x+y)+z=x+(y+z)$ for all $x,y,z \in H$,
\item[(CH2)] (commutativity) $x+y=y+x$ for all $x,y\in H$,
\item[(CH3)] (inverse element) for every $x\in H$ there exists a unique $x'\in H$ such that $0\in x+x'$ (the element $x'$ will be denoted by $-x$),
\item[(CH4)] (reversibility axiom) $z\in x+y$ implies $x\in z-y:=z+(-y)$ for all $x,y,z\in H$.
\end{itemize}
\end{definition}

\begin{remark}
We call $0$ the neutral element for $+$. In the literature, canonical hypergroups often require explicitly that $x+0 = \lbrace x \rbrace$ for all $x \in H$. This is a consequence of the Axioms (CH3) and (CH4): suppose that $y \in x+0$ for some $x,y \in H$. Then $0 \in y-x$ by (CH4). Now $y = x$ follows from the uniqueness required in (CH3). 
\end{remark}

\begin{remark}
A canonical hypergroup is a hypergroup in the sense of Marty. Fix $a \in H$ and take $x \in H+a$. Then there exist $h \in H$ such that $x \in h+a \subseteq H$, showing that $H+a \subseteq H$. For the other inclusion, take $x \in H$, then 
\[
x \in x+0 \subseteq x+(a-a) = (x-a) + a,
\]
so there exists $h \in x-a \subseteq H$ such that $x \in h+a \subseteq H+a$.
\end{remark}

\begin{definition}\label{Krasner'shyperrings}
A \emph{hyperfield} is a structure $\mathcal{H}=(H,+,\cdot,0,1)$ which satisfies the following axioms:
\begin{itemize}
\item[(HF1)] $(H,+,0)$ is a canonical hypergroup,
\item[(HF2)] $(H,\ \cdot, 0,1 \ )$ is a commutative pseudo-group, i.e., $H\setminus\{0\}$ is an abelian group and $0$ is an absorbing element: $x\cdot0= 0$ for all $x\in H$,
\item[(HF3)] (distributivity of $\cdot$ over $+$) $x\cdot(y+z)=x\cdot y+x\cdot z$ for all $x,y,z\in H$.
\end{itemize}
In (HF3), for $x\in H$ and $A\subseteq H$, we have set 
\[ 
x\cdot A:=\{x\cdot a\mid a\in A\}.
\]
\end{definition}

\begin{remark} \label{RemarkDoubleDistributivity}
As was shown by Viro in \cite[Section 4.4]{Vir10}, the double distributivity law, i.e.,
\[
(a+b)(c+d) = ac+ad+bc+bd,
\]
does not hold in general in hyperfields. However, we always have the following inclusion:
\[
(a+b)(c+d) \subseteq ac+ad+bc+bd.
\] 
\end{remark}

\begin{example}
Let $\mathbb{K}:=\{0,1\}$. Define an hyperoperation $+$ on $\mathbb{K}$ by setting $0$ as the neutral element and $1+1:=\{0,1\}$. Then, equipped with the standard multiplication $\cdot$, the structure $(\mathbb{K}, +,\cdot, 0,1)$ is a hyperfield.\par
Let $\mathbb{S}:=\{-1,0,1\}$. Define a hyperoperation $+$ on $\mathbb{S}$ by setting $0$ as the neutral element and $1+1:=\{1\}$, $-1-1:=\{-1\}$ and $1-1:=\{-1,0,1\}$. Then, equipped with the standard multiplication $\cdot$, the structure $(\mathbb{S}, +,\cdot, 0,1)$ is a hyperfield.
\end{example}

Some hyperfields (such as $\mathbb{K}$ and $\mathbb{S}$) can be very close to fields, in a sense which is captured by the following definition.

\begin{definition}[\cite{BS21}]
A hyperfield $\mathcal{H}$ is called \emph{stringent} if for all $x,y \in \mathcal{H}$ we have that $x+y$ is a singleton, unless $0\in x+y$.    
\end{definition}

Next, we introduce the notion of valued hyperfield in the sense of Krasner. We will see that when these hyperfields are stringent, their structure can be described in several interesting ways (cf.\ Fact~\ref{biinterpretability}).

\begin{definition}[\cite{Kra57}]\label{DefinitionValuationHyperfields}
     Let $\Gamma= (\Gamma, +,0,<)$ be an ordered abelian group and let $\rho$ be a non-empty initial segment of $\Gamma_{\geq 0}$. A \emph{valuation map (of norm $\rho$)} on a hyperfield $\mathcal{H}=(H,+,\cdot,0,1)$ is a surjective function $\val: H \rightarrow \Gamma\cup\{\infty\}$ satisfying for all $x,y\in H$

\begin{itemize}
    \item[(V0)] $\val(x)= 0$ if and only if $x=0$.
    \item[(V1)] $\val(xy)=\val(x)+\val(y)$,
    \item[(V2)] (Ultrametric inequality) $\val(z) \geq \min(\val(x),\val(y))$ for all $z\in x+ y$,
    \item[(V3)] $0 \notin x+ y \Longrightarrow  \vert \val(x+ y) \vert =1$,
    \item[(V4)] (Norm axiom) For all $z,z' \in H$ we have that $z\in x+ y$ implies that 
    \[
    z'\in x+ y \iff \forall w\in z-z'~:~\val(w)>\rho+ \min\{\val(x),\val(y)\}.
    \]
\end{itemize}
As usual, $\infty$ is understood as a symbol such that $\infty>\Gamma$ and $\gamma+\infty=\infty+\gamma=\infty$ for all $\gamma\in\Gamma$. A \emph{valued hyperfield} is a hyperfield equipped with a valuation map. If $\mathcal{H}$ is a valued hyperfield, then we denote its norm by $\mathcal{N}(\mathcal{H})$.
\end{definition}

\begin{notation}
Note that in the definition above we have used the symbol $+$ to denote both the hyperoperation of $\mathcal{H}$ as well as the operation of $\Gamma$ as it is customary in classical valuation theory on fields. If there is risk of confusion, then we will prefer the symbol $\boxplus$ to denote hyperoperations.    
\end{notation}

\begin{remark}
    In the literature (see for instance \cite{KLS22}) one may find a weaker notion of valuation on hyperfields. This is obtained by asking for axioms (V0), (V1) and (V2) only.
\end{remark}
The next lemma contains some properties of valuations on hyperfields which are analogous to those of classical valuations on fields. Note that these follow from axioms (V0), (V1) and (V2) only.
\begin{lemma}[Lemma 4.5 in \cite{KLS22}]\label{propertiesofvaluations}
Let $\val:\mathcal{H}\to\Gamma\cup\{\infty\}$ be a valuation on a hyperfield $\mathcal{H}$. Then:
\begin{enumerate}
 \item[$(i)$] $\val(1)=\val(-1) = 0$,
 \item[$(ii)$] $\val(-x)=\val(x)$ for all $x\in \mathcal{H}$,
 \item[$(iii)$] $\val(x^{-1})=-\val(x)$ for all $x\in \mathcal{H}$,
 \item[$(iv)$]  For all $x,y\in\mathcal{H}$, if $\val(x)\neq \val(y)$, then 
 $\val(z) = \min\{\val(x), \val(y)\}$ for every $z \in x+y$.
\end{enumerate}
\end{lemma}

\begin{remark}\label{KrasnerValuedUltrametric}
Let us explain the role that axioms (V3) and (V4) have by comparing valued hyperfields with valued fields. When a valued field with valuation $\val$ is given, then $d(x,y):=\val(x-y)$ always defines an ultrametric on the field. As in a hyperfield $x-y$ denotes a set, Axiom (V3) is intended to ensure that the analogous function $d: (x,y) \mapsto \val(z)$ for some $z\in x-y$ gives a well-defined ultrametric (induced by the valuation). With respect to this ultrametric, Axiom (V4) postulates that $x+y$ must be an ultrametric ball of radius $\rho+\min\{\val(x),\val(y)\}$, where $\rho$ is the norm of the valuation. This provides some control on the subsets that can be obtained as the sum of two elements of the hyperfield. 
\end{remark}

Let us observe a useful property of valued hyperfields which can be deduced from Axiom (V4).

\begin{lemma}\label{RemarkSumIsSingleton}
Let $\mathcal{H}$ be a valued hyperfield with valuation $\val$ and $x,y\in\mathcal{H}$. If $x\in x+y$, then $x+y=\{x\}$.
\end{lemma}

\begin{proof}
Indeed, by reversibility, if $x\in x+y$, then $y\in x-x$ and so $\val(y)>\mathcal{N}(\mathcal{H})+\val(x)\geq \val(x)$ by (V4). Therefore, if $z\in x+y$, then $\val(x)=\val(z)$ by Lemma \ref{propertiesofvaluations} $(iv)$. We have $\val(y-0)>\mathcal{N}(\mathcal{H})+\min\{\val(z),\val(-x)\}$ and since $y\in z-x$ (by reversibility), by (V4) again, we obtain that $0\in z-x$. This implies that $x=z$ holds and proves the assertion of the lemma.
\end{proof}

\begin{remark}
We remark that, for instance, the hyperfield $\mathbb{S}$ does not satisfy the above property. Indeed, $1\in1-1$ but $1-1\neq\{1\}$ in $\mathbb{S}$.
\end{remark}

We show in the next proposition that Axiom (V2) and (V4) generalise for sums of elements of arbitrary length:
\begin{proposition}\label{PropositionSumOfNElementsIsABall}
Let $\mathcal{H}$ be a valued hyperfield with valuation $\val$.
\begin{itemize}
    \item[$(i)$] For all $y\in x_1+\cdots+x_n$ we have that $\val(y) \geq \min\{\val(x_i)\mid i=1,\ldots,n\}$ for all $n\in\mathbb{N}_{\geq 2}$. 
\item[$(ii)$] Let $n\in\mathbb{N}_{\geq2}$. For all $x_1,\ldots,x_n,y,y' \in H$ we have that $y\in x_1+\ldots+x_n$ implies that $y'\in x_1+\ldots+x_n$ if and only if $\val(w)>\rho+ \min\{\val(x_i)\mid i=1,\ldots,n\}$ for all $w\in y-y'$. 
\end{itemize}
\end{proposition}
\begin{proof}\
    \begin{itemize}
        \item[$(i)$] We prove this by induction on $n$. The base step is given by axiom (V2). Let $n>2$ and take $y\in x_1+\cdots+x_n$. We can assume without loss of generality that $\min\{\val(x_i)\mid i=1,\ldots,n\}=\val(x_1)$. Let $z\in x_2+\cdots+x_n$ be such that $y\in x_1+z$. By the induction hypothesis we obtain that $\val(z)\geq\min\{\val(x_i)\mid i=2,\ldots,n\} \geq \val(x_1)$. Thus, by Axiom (V2), we obtain that $\val(y)\geq\min\{\val(x_1),\val(z)\}\geq\val(x_1)$, as we wished to show.
        \item[$(ii)$] We proceed by induction on $n\in \mathbb{N}_{\geq2}$. The base step is axiom (V4). For the induction step, let $n>2$ and assume without loss of generality that $\min\{\val(x_i)\mid i=1,\ldots,n\}=\val(x_1)$. Take $y\in x_1+ \dots + x_n$. Let $z\in x_2+\cdots+x_n$ be such that $y\in x_1+z$. Now, if $y'\in \mathcal{H}$ is taken such that $\val(w)>\rho+\val(x_1)$ for all $w\in y-y'$, then for all $z'\in y'-x_1$, we have that $z'-z\subseteq (y'-y)+(x_1-x_1)$ and by Axiom (V2) we obtain that $\val(u) > \rho +\min\{\val(x_i)\mid i=1,\ldots,n\}$ for all $u\in z-z'$. By the induction hypothesis, it follows that $z'\in x_2+\cdots+x_n$ and hence $y'\in x_1+\cdots+x_n$ by reversibility.
Conversely, if $y'\in x_1+ \cdots + x_n$, then there exists $z'\in x_2+\cdots + x_n$ such that $y'\in x_1+z'$. By the induction hypothesis $\val(u) > \rho +\min\{\val(x_i)\mid i=2,\ldots,n\}$ for all $u\in z-z'$. On the other hand, we have that $y-y'\subseteq (x_1-x_1)+(z-z')$, so, using (V2) and (V4) again, we may conclude that $\val(w)>\rho+\val(x_1)$ for all $w\in y-y'$.\qedhere
    \end{itemize}
\end{proof}

\subsection{Valued hyperfields associated to valued fields}
    
We will recall now how we can canonically associate valued hyperfields to a valued field.\par
First, consider a field $\mathcal{F}=(F,+,\cdot,0,1)$ and let $T$ be a subgroup of $\mathcal{F}^\times$. For $x\in \mathcal{F}^\times$ we denote by $[x]_T$ the coset $xT\in \mathcal{F}^\times/T$. Further, let $[0]_T$ denote the singleton containing only $0\in \mathcal{F}$. Then the \emph{factor hyperfield} of $\mathcal{F}$ modulo $T$ is the set $\mathcal{F}_T:=\mathcal{F}^\times/T\cup\{[0]_T\}$ with the hyperoperation (for which we again use the symbol $+$) defined as:
\[
[x]_T+[y]_T:=\{[x+yt]_T\mid t\in T\}
\]
and with the obvious multiplication
\[
[x]_T[y]_T:=[xy]_T~.
\]
This construction was shown to yield always a hyperfield by Krasner in \cite{Kra83}. In that paper, Krasner also conjectured that all hyperfields are of this form. This was later shown to be false by Massouros in \cite{Mas85}.
Before we go on, we need to generalise Krasner's construction, starting from a hyperfield.

\begin{proposition}
    Let $\mathcal{H}=(H,+,\cdot,0,1)$ be a hyperfield, and $T$ a subgroup of $\mathcal{H}^\times$. As before denote
$[x]_T$ the coset $xT\in \mathcal{H}^\times/T$ and by $[0]_T = \{0\}$. Define
\[
[x]_T+[y]_T:=\{[z]_T\mid  \exists t\in T~:~z\in x+yt\}.
\]
Then $\mathcal{H}_T:= (\mathcal{H}^\times/T\cup \{[0]_T\}, +, \cdot, [0]_T, [1]_T )$ is a hyperfield (where $\cdot$ is the obvious multiplication).  
\end{proposition}

\begin{proof}

We begin by showing that axiom (CH1) holds in $\mathcal{H}_T$. Let $[x]_T,[y]_T$ and $[z]_T$ be arbitrary elements of $\mathcal{H}_T$ (and therefore $x,y,z\in\mathcal{H}$). We have that
    \begin{equation}\label{ass1}
         ([x]_T+[y]_T)+[z]_T =\bigcup_{[u]_T\in[x]_T+[y]_T}[u]_T+[z]_T=\bigcup_{[u]_T\in[x]_T+[y]_T}\{[v]_T\mid\exists t\in T~:~v\in u+zt\}
    \end{equation}
and
 \begin{equation}\label{ass2}
         [x]_T+([y]_T+[z]_T) =\bigcup_{[u]_T\in[y]_T+[z]_T}[x]_T+[u]_T=\bigcup_{[u]_T\in[y]_T+[z]_T}\{[v]_T\mid\exists t\in T~:~v\in x+ut\}.
    \end{equation}
We must show that these two sets are equal and to do so we will show two inclusions.\par
In order to show that \eqref{ass1} is contained in \eqref{ass2}, pick $[v]_T$ such that there exists $t\in T$ with $v\in u+zt$ for some $u$ such that there exists $s\in T$ with $u\in x+ys$. Then $v\in (x+ys)+zt$. By distributivity and associativity in $\mathcal{H}$ we may conclude that $vs^{-1}\in xs^{-1}+(y+zts^{-1})$. Therefore, there exists $w\in y+zts^{-1}$ such that $vs^{-1}\in xs^{-1}+w$ and hence $v\in x+ws$, where we again used distributivity in $\mathcal{H}$. This shows that $[v]_T$ is an element of \eqref{ass2}. For the converse inclusion, pick $[v]_T$ such that there exists $t\in T$ with $v\in x+ut$ for some $u$ such that there exists $s\in T$ with $u\in y+zs$. Then using distributivity and associativity in $\mathcal{H}$ we obtain that
\[
v\in x+(y+zs)t=(x+yt)+zst.
\]
Therefore, there exists $w\in x+yt$ such that $v\in w+zst$ and this shows that $[v]_T$ is an element of \eqref{ass1}.\par    
Next, we show that axiom (CH2) holds in $\mathcal{H}_T$. Let $[x]_T$ and $[y]_T$ be arbitrary elements of $\mathcal{H}_T$. Pick $[z]_T\in [x]_T+[y]_T$ so that there exists $t\in T$ such that $z\in x+yt$. By distributivity and commutativity in $\mathcal{H}$ we obtain that $zt^{-1}\in y+xt^{-1}$ and therefore $[z]_T=[zt^{-1}]_T\in [y]_T+[x]_T$ follows immediately. This suffices to show that the hyperoperation of $\mathcal{H}_T$ is commutative.\par
We now show that axiom (CH3) holds in $\mathcal{H}_T$. Let $[x]_T$ be an arbitrary element of $\mathcal{H}_T$ and assume that $[0]_T\in[x]_T+[y]_T$ for some $[y]_T\in \mathcal{H}_T$. This means that there exists $t\in T$ such that $0\in x+yt$ and, since axiom (CH3) holds in $\mathcal{H}$, we may conclude that $yt=-x$ and so $[y]_T=[-x]_T$. Thus, $[-x]_T$ is the unique element $[y]_T\in \mathcal{H}_T$ such that $[0]_T\in [x]_T+[y]_T$. In symbols, we have shown that $-[x]_T=[-x]_T$.\par
Finally, we prove that axiom (CH4) holds in $\mathcal{H}_T$. Let $[x]_T$ and $[y]_T$ be arbitrary elements of $\mathcal{H}_T$ and assume that $[z]_T\in[x]_T+[y]_T$. We have to show that $[x]_T\in[z]_T-[y]_T$. By assumption there exists $t\in T$ such that $z\in x+yt$. By reversibility in $\mathcal{H}$, we obtain that $x\in z-yt=z+(-y)t$. Hence, $[x]_T\in[z]_T-[y]_T$ since $-[y]_T=[-y]_T$ as we have proved above.\par

By definition, $(\mathcal{H}_T\setminus\{[0]_T\},\cdot,[1]_T)$ is a group so that axiom (HF2) holds in $\mathcal{H}_T$. It remains to show that axiom (HF3) holds in $\mathcal{H}_T$. For this let $[x]_T$, $[y]_T$ and $[z]_T$ be arbitrary elements of $\mathcal{H}_T$. Using distributivity in $\mathcal{H}$, we can make the following straightforward computation:
    \begin{align*}
         [z]_T\cdot ([x]_T+[y]_T) & =\{[zv]_T\mid  \exists t\in T~:~v\in x+yt\} \\
         & = \{[w]_T\mid \exists t\in T~:~ w\in z(x+yt)\} \\
         & = \{[w]_T\mid \exists t\in T~:~ w\in zx+zyt\}\\
         & =[zx]_T+[zy]_T\\
         & =[z]_T[x]_T+[z]_T[y]_T.
    \end{align*}
We have proved that $\mathcal{H}_T$ is a hyperfield.
\end{proof}

\begin{remark}\label{RemarkClassofSumIsAlwaysInSumClass}
    For all $x,y\in \mathcal{H}$, since $1\in T$, by definition we have that $[z]_T \in [x]_T+[y]_T$ for all $z\in x+y$. In particular, if $\mathcal{H}$ is a field, then $[x+y]_T\in[x]_T+[y]_T$ always.
\end{remark}

Consider now a valued field $\mathcal{K}=(K,+,\cdot,0,1,\val)$ with value group $\Gamma$.
If $\rho$ is an initial segment of $\Gamma_{\geq 0}$, then we denote by $\mathfrak{m}_\rho$ the ideal of the valuation ring of $\mathcal{K}$ consisting of all elements of value greater than $\rho$. Observe that the elements of $T:=1+\mathfrak{m}_\rho$ have all value $0$ and thus the valuation of $\mathcal{K}$ is well-defined as a map on $\mathcal{K}_T$. Moreover, it is not difficult to show that this map is a valuation on $\mathcal{K}_T$ of norm $\rho$.

\begin{notation}
We denote the valued hyperfield $\mathcal{K}_{1+\mathfrak{m}_\rho}$ by $\mathcal{H}_\rho(\mathcal{K})$ and its valuation by $\val_\rho$. In this case, we write $[x]_\rho$ in place of $[x]_{1+\mathfrak{m}_\rho}$ and $\theta_\rho: \mathcal{K} \rightarrow \mathcal{H}_\rho(\mathcal{K})$ for the canonical epimorphism $x\mapsto[x]_\rho$. 
\end{notation}

We also need to generalise this construction, starting with a valued hyperfield. 
\begin{proposition}\label{PropositionQuotientHyperfields}
    Let $\mathcal{H}:=(H,+,\cdot,0,1,\val)$ be a valued hyperfield with value group $\Gamma$. For any initial segment $\rho \subseteq \mathcal{N}(\mathcal{H})$, let $\mathfrak{m}_\rho$ be the set of all the elements of $\mathcal{H}$ with value strictly bigger than $\rho$. Then 
    \begin{itemize}
        \item[$(i)$] $T_\rho:=1+\mathfrak{m}_\rho$ is a subgroup of $\mathcal{H}^\times$.
        \item[$(ii)$] The valuation on $\mathcal{H}$ induces a valuation $\val_\rho:\mathcal{H}_\rho \rightarrow \Gamma$ on the factor hyperfield $\mathcal{H}_\rho := \mathcal{H}_{T_\rho}$.
        \item[$(iii)$] The valuation $\val_\rho$ is of norm $\rho$.
    \end{itemize}
\end{proposition}
    Notice that if $\rho$ is an initial segment of $\Gamma_{\geq 0}$ containing $\mathcal{N}$, then by Axiom (V4) we have $1+\mathfrak{m}_\rho=\{1\}$ and so $\mathcal{H}_\rho = \mathcal{H}$. 

\begin{proof}
    As for valued fields, the fact that $1+\mathfrak{m}_\rho$ is a subgroup of $\mathcal{H}^\times$ follows easily from the axioms of valuations: if $x,y\in1+\mathfrak{m}_\rho$, then we have
    $xy\in(1+t)(1+s)$ for some $t,s\in\mathfrak{m}_\rho$, and by the observation made in Remark \ref{RemarkDoubleDistributivity} we obtain that $xy\in 1+(t+s+ts)$. Then there is $z\in t+s+ts$ such that $xy\in 1+z$ where $z\in\mathfrak{m}_\rho$ by Proposition \ref{PropositionSumOfNElementsIsABall} $(i)$.
    Let us now show that
    \begin{align*}
    \val_\rho: \mathcal{H}_\rho&\to\Gamma\cup\{\infty\}\\
    [x]_\rho &\mapsto \val(x)
    \end{align*}
    defines a valuation of norm $\rho$ on $\mathcal{H}_\rho$. By Lemma \ref{propertiesofvaluations} $(iv)$, the elements of $1+\mathfrak{m}_\rho$ have value $0$, so $\val_\rho$ is a well-defined and surjective map. We further have that
    \begin{itemize}
        \item[(V0)] $\val_\rho([x]_\rho)= \infty  \ \Longleftrightarrow \val(x)=\infty \Longleftrightarrow x= 0 \Longleftrightarrow [x]_\rho=[0]_\rho$.
        \item[(V1)]
        $\val_\rho([x]_\rho  [y]_\rho)=\val_\rho([xy]_\rho )=\val(xy)=\val(x)+\val(y)=\val_\rho([x]_\rho)+\val_\rho([y]_\rho)$.
        \item[(V2)] For all $[x]_\rho,[y]_\rho\in\mathcal{H}_\rho$, if $[z]_\rho \in [x]_\rho+[y]_\rho$, then by definition there exists $t\in 1+\mathfrak{m}_\rho$ such that $z\in x+yt$ and since axiom (V2) holds for $\val$, we have that 
        \[
        \val_\rho([z]_\rho)=\val(z)\geq\min\{\val(x),\val(yt)\}=\min\{\val(x),\val(y)\}=\min\{\val_\rho([x]_\rho),\val_\rho([y]_\rho)\}.
        \]
        \item[(V3)] Take $[x]_\rho,[y]_\rho\in\mathcal{H}_\rho$ and assume without loss of generality that $\val_\rho([y]_\rho)\leq\val_\rho([x]_\rho)$. If $| \val([x]_\rho+[y]_\rho)|>1$, then there exist $t\in 1+\mathfrak{m}_\rho$, $z\in x+y$ and $z'\in x+yt$ such that $\val(z')>\val(z)$. We claim that $z'\in x+y$. Indeed, if $w\in z-z'$, then $w\in (x-x)+(y-ty)$, so there are $a\in x-x$ and $b\in (1-t)y$ such that $w\in a+b$. By axiom (V4) applied in $\mathcal{H}$ and since $t\in 1+\mathfrak{m}_\rho$, we deduce that $\val(a)>\rho+\val(x)\geq\rho+\val(y)$ and that $\val(b)>\rho+\val(y)$. Therefore, by axiom (V2) applied in $\mathcal{H}$, we obtain that
        \[
        \val(w)\geq\min\{\val(a),\val(b)\}>\rho+\val(y)=\rho+\min\{\val(x),\val(y)\}.
        \]
        Hence, axiom (V4) applied in $\mathcal{H}$ now implies that $z'\in x+y$. Since $z$ and $z'$ have distinct values, by axiom (V3) applied in $\mathcal{H}$ we conclude that $x=-y$ must hold and thus $[x]_\rho=-[y]_\rho$. This proves that (V3) holds in $\mathcal{H}_\rho$. 
        \item[(V4)] Take $[x]_\rho,[y]_\rho\in\mathcal{H}_\rho$ and assume without loss of generality that $\val_\rho([y]_\rho)\leq\val_\rho([x]_\rho)$. If $[z]_\rho \in [x]_\rho + [y]_\rho$, then there exists $t\in 1+\mathfrak{m}_\rho$ such that $z \in x+ty$. Take $[z']_\rho\in[x]_\rho+[y]_\rho$ so that there exists $t'\in1+\mathfrak{m}_\rho$ such that $z'\in x+t'y$. Fix $[w]_\rho\in[z]_\rho-[z']_\rho$. There is $u\in1+\mathfrak{m}_\rho$ such that 
        \[
        w\in z-uz'\subseteq (x+ty)-(x+t'y)u=(x-xu)+(yt-yt'u)
        \]
        and hence we may find $a\in (1-u)x$ and $b\in (t-t'u)y$ such that $w\in a+b$. By axiom (V4) applied in $\mathcal{H}$ and since the elements of $1-u$ and $t-t'u$ have value $>\rho$, we obtain that $\val(a)>\rho+\val(x)\geq\rho+\val(y)$ and that $\val(b)>\rho+\val(y)$. Therefore, axiom (V2) applied in $\mathcal{H}$ yields
        \[
        \val(w)\geq\min\{\val(a),\val(b)\}>\rho+\val(y)=\rho+\min\{\val(x),\val(y)\}.
        \]
        This shows that 
        \[
        \val_\rho([w]_\rho)>\rho+\min\{\val_\rho([x]_\rho),\val_\rho([y]_\rho)\}
        \]
        for all $[w]_\rho\in[z]_\rho-[z']_\rho$.\par
        Conversely, assume that $\val_\rho([w]_\rho)>\rho+\val_\rho([y]_\rho)$ for all $[w]_\rho\in[z]_\rho-[z']_\rho$, where $[z']_\rho\in\mathcal{H}_\rho$. We have to show that $[z']_\rho\in[x]_\rho+[y]_\rho$. If $[z]_\rho=[z']_\rho$, then there is nothing to prove. Otherwise, take any $w\in z-z'$. By assumption we have that $\val(w)>\rho+\val(y)$ and hence $wy^{-1}\in\mathfrak{m}_\rho$. Let us now pick $a\in\mathfrak{m}_\rho$ such that $t\in 1+a$. By reversibility in $\mathcal{H}$ we obtain that $z'\in z-w\subseteq x+yt-w\subseteq x+y+y(a-wy^{-1})$. Hence, there exists $b\in a-wy^{-1}$ such that $z'\in x+y(1+b)$. Now, by axiom (V2) applied in $\mathcal{H}$, we have that $\val(b)\geq\min\{\val(a),\val(wy^{-1})\}>\rho$. Thus, there is $u\in 1+b\subseteq1+\mathfrak{m}_\rho$ such that $z'\in x+yu$, i.e., $[z']_\rho\in[x]_\rho+[y]_\rho$.
    \end{itemize}
We have proved that $\val_\rho$ is a valuation of norm $\rho$ on $\mathcal{H}_\rho$.    
\end{proof}

\subsection{The inverse limit of some isometric systems of valued hyperfields}    
    
Let us recall another notion due to Krasner.
    
    \begin{definition}[\cite{Kra57}]\label{DefinitionIsometry}
    Let $\mathcal{H}$ and $\mathcal{H}'$ be valued hyperfields with the same value group $\Gamma$. An \emph{isometric homomorphism} from $\mathcal{H}$ to $\mathcal{H}'$ is a map $\theta:\mathcal{H} \rightarrow \mathcal{H}'$ satisfying:
    \begin{enumerate}
        \item[(IH1)] $\theta(x y)=\theta(x)  \theta(y)$ for all $x,y \in \mathcal{H}$,
        \item[(IH2)]  $\theta^{-1}(x'+y')= \theta^{-1}(x') + \theta^{-1}(y')$ for all $x',y' \in \theta(\mathcal{H})$,
        \item[(IH3)] $\val(x)=\val'(\theta(x))$ for all $x \in \mathcal{H}$.
    \end{enumerate}
\end{definition}
\begin{remark}
Notice that Krasner originally considered only valuations of rank one. This means that the value groups were all subgroups of $\mathbb{R}$. To make sense of axiom (IH3) in the general case, one can assume (as we do) that the valued hyperfields have the same value group. More generally, one may ask that these value groups both \lq\lq live\rq\rq\ in a (possibly bigger) ordered abelian group (which in the rank one case is $\mathbb{R}$). Here, we will not elaborate more on this aspect.\par
\end{remark}

The above definition is inspired by the following situation. 

\begin{example}\label{canonicalisometricsystem}
Let $\mathcal{K}$ be a valued field with value group $\Gamma$ and let $\rho'\subseteq\rho$ be initial segments of $\Gamma_{\geq0}$. It is not difficult to check that the canonical epimorphism $\theta_{\rho'}:\mathcal{K}\to\mathcal{H}_{\rho'}(\mathcal{K})$ induces on $\mathcal{H}_\rho(\mathcal{K})$ a surjective isometric homomorphism $\theta_{\rho,\rho'}:\mathcal{H}_{\rho}(\mathcal{K})\to\mathcal{H}_{\rho'}(\mathcal{K})$.
\end{example}

Let us observe that in some situations, we can replace axiom (IH2) with the following:

\begin{itemize}
    \item[(IH2')] $\theta(x+y)\subseteq\theta(x)+\theta(y)$ for all $x,y\in \mathcal{H}$.
\end{itemize}

In the next lemma, we describe such a situation.

\begin{lemma}\label{PropositionIH2EquivalentIH2'}
Let $\mathcal{H}$ and $\mathcal{H}'$ be hyperfields. \begin{itemize}
    \item[$(i)$] If a map $\theta: \mathcal{H}\to \mathcal{H}'$ satisfies (IH2) then it also satisfies (IH2'). 
    \item[$(ii)$] Assume that $\mathcal{H}'$ is a valued hyperfield with valuation $\val'$. If $\theta: \mathcal{H} \to \mathcal{H}'$ is a map satisfying $\theta(0)=0$ and (IH2'), then it satisfies (IH2).
    \end{itemize}
\end{lemma}

\begin{proof}
We show $(i)$. Let $x,y\in \mathcal{H}$ and pick $z\in x+y$. We have to show that $\theta(z)\in \theta(x)+\theta(y)$. If $x'=\theta(x)$ and $y'=\theta(y)$, then $z\in x+y$ implies $z\in\theta^{-1}(x')+\theta^{-1}(y')$. By (IH2), $z\in\theta^{-1}(x'+y')$ which means $\theta(z)\in x'+y'=\theta(x)+\theta(y)$.\par 
Now we show $(ii)$. Using (IH2'), we obtain that $0=\theta(0) \in \theta(y-y) \subseteq \theta(y)+\theta(-y)$. Then, by the uniqueness required in axiom (CH3), we must have $-\theta(y)=\theta(-y)$. To show that (IH2) holds, take $x',y'\in\theta(\mathcal{H})$. If $z\in\theta^{-1}(x')+\theta^{-1}(y')$, then there are $x\in\theta^{-1}(x')$ and $y\in\theta^{-1}(y')$ such that $z\in x+y$. It follows by (IH2') that 
\[
\theta(z)\in\theta(x+y)\subseteq\theta(x)+\theta(y)=x'+y'
\] 
and hence $z\in \theta^{-1}(x'+y')$. For the other inclusion, take $z\in\theta^{-1}(x'+y')$. Assume without loss of generality that $\val'(x')\leq \val'(y')$. If $y\in\theta^{-1}(y')$, then by (IH2') we have $\theta(z-y)\subseteq\theta(z)-\theta(y)\subseteq (x'+y')-y'=x'+(y'-y')$. If $a\in y'-y'$, then $\val'(a)>\mathcal{N}(\mathcal{H}')+\val'(y')\geq\mathcal{N}(\mathcal{H}')+\val'(x')$. It follows that $a\in x'-x'$ and $x'\in x'+a=\{x'\}$ by Remark \ref{RemarkSumIsSingleton}. We have shown that for all $a\in y'-y'$ we have that $x'+a=\{x'\}$. Therefore, $x'+(y'-y')=\{x'\}$. Thus, $z-y\subseteq\theta^{-1}(x')$ for all $y\in\theta^{-1}(y')$. By reversibility, this proves that $z\in\theta^{-1}(x')+\theta^{-1}(y')$ as we wished to show.
\end{proof}

\begin{remark}
In the literature, a homomorphism of (not necessarily valued) hyperfields $\theta:\mathcal{H}\to\mathcal{H}'$ is commonly defined to be a map satisfying $\theta(0)=0$, (IH1) and (IH2'). Part $(i)$ of the above lemma shows that any isometric homomorphism is in fact a homomorphism of the underlying hyperfields.
\end{remark}

\begin{remark}
Note that $\theta(0)=0$ is always satisfied by an isometric homomorphism as it follows from (V0) and (IH3). This means that when considering isometric homomorphisms between valued hyperfields properties (IH2) and (IH2') are equivalent.
\end{remark}

\begin{remark}
If $\theta:\mathcal{H}\to \mathcal{H}'$ is an isometric homomorphism between two valued hyperfields $\mathcal{H}$ and $\mathcal{H}'$, then $\mathcal{N}(\mathcal{H}) \geq\mathcal{N}(\mathcal{H}')$. This follows from Axiom (V4) because $\theta(1_\mathcal{H}-1_\mathcal{H})\subseteq 1_{\mathcal{H}'}-1_{\mathcal{H}'}$ by (IH2'). Indeed, the above inclusion means that for all $x\in \mathcal{H}$, if $\val(x)>\mathcal{N}(\mathcal{H})$, then by (IH3) $\val(x)=\val'(\theta(x))>\mathcal{N}(\mathcal{H}')$.
\end{remark}

\begin{lemma}\label{LemmaIsometryIndicesIsomophism}
    Let $\mathcal{H}$ and $\mathcal{H}'$ be valued hyperfields, and $\theta$ be an isometric homomorphism from $\mathcal{H}$ to $\mathcal{H}'$. Then $\theta$ induces an isomorphism of valued hyperfields 
    \[
    \theta': \mathcal{H}_{\mathcal{N}'} \rightarrow \mathcal{H}', 
    \]
    where $\mathcal{N}'$ is the norm of $\mathcal{H}'$ and $\mathcal{H}_{\mathcal{N}'}$ is the hyperfield of norm $\mathcal{N}'$ obtained from $\mathcal{H}$ as in Proposition \ref{PropositionQuotientHyperfields}. 
\end{lemma}
\begin{proof}
    We remark that $T:=\theta^{-1}(1)$ is a multiplicative subgroup of $\mathcal{H}$. The isometry $\theta$ induces an isomorphism 
    \[\theta': \mathcal{H}_{T} \rightarrow \mathcal{H}', [a]_T \mapsto \theta(a) \]
    Indeed, it is well defined and it is a isomorphism of group by (IH1) and preserves the valuation by (IH3). To prove the isomorphism of valued hyperfields, it remains to show that it preserves the additive structure. Notice that $\theta^{-1}(\theta(a))=[a]_T$ for all $a\in \mathcal{H}$. Then we have:
    \begin{align*}
        [c]_T\in  [a]_T+[b]_T &\Leftrightarrow c \in  \theta^{-1}(\theta(a))+\theta^{-1}(\theta(b))\\
            & \Leftrightarrow  c \in \theta^{-1}(\theta(a)+\theta(b)) \qquad \text{by (IH2)},\\
            & \Leftrightarrow  \theta(c) \in  \theta(a)+\theta(b).
    \end{align*}
    Since $\mathcal{H}$ is a valued hyperfield, it remains to show that $T=1+\mathfrak{m}_{\mathcal{N}'}$. Let $x\in \theta^{-1}(1)$, and  $y\in x-1$. By Axiom (CH4) it is enough to show that $\val(y)>\mathcal{N}'$. By (IH2), $\theta(y)\in \theta(x)-1=1-1$. Thus, by $\text{(V4)}_{\mathcal{H'}}$, we have  $\theta(y)>\mathcal{N}'$ and thus by (IH3) $\val(y)>\mathcal{N}'$. Conversely, if $x\in1+y$ with $y\in \mathfrak{m}_{\mathcal{N}'}$, then $\theta(x)\in 1+\theta(y)$ with $\val(\theta(y))=\val(y)>\mathcal{N}'$. Then by Axiom (V4) $\theta(y) \in 1-1$. By (CH4) $1\in 1+\theta(y)$ and then, Lemma \ref{RemarkSumIsSingleton} yields  $1+\theta(y)=\{1\}$. Thus, we have $\theta(x)=1$ as desired. 
\end{proof}

Example \ref{canonicalisometricsystem} above inspires the following definition.

\begin{definition}
    Let $(I,<)$ be a total order. We say that a family of valued hyperfields $\{\mathcal{H}_i\mid i\in I\}$ with the same value group $\Gamma$ is an \emph{isometric system} if for all $i<j$ there is a surjective isometric homomorphism $\theta_{j,i}:\mathcal{H}_j\to\mathcal{H}_i$ and $\theta_{j,i}\circ\theta_{k,j}=\theta_{k,i}$ for all $i<j<k$.
\end{definition}

In the following proposition, we describe the structure of the inverse limit of isometric systems of valued hyperfields, satisfying a certain assumption. Again, the rank one case was already treated by Krasner in \cite[Paragraph 4]{Kra57}. We adapt and extend his results to fit the context of non-archimedian value groups.

\begin{proposition}\label{PropositionInverseLimitHyperfields}
   Let $(I,<)$ be a total order and  $(\{\mathcal{H}_i\mid i\in I\},\{\theta_{j,i}\mid i<j\})$ an isometric system of valued hyperfields with value group $\Gamma$. If the sequence of norms $(\mathcal{N}(\mathcal{H}_i))_i$ is cofinal in a convex subgroup of $\Gamma$ (equivalently, for all $i$ there is $j$ such that $\mathcal{N}(\mathcal{H}_{j})\geq 2 \mathcal{N}(\mathcal{H}_{i})$), then
    
    \begin{enumerate}
        \item[(1)] The inverse limit $\mathcal{H}:= \varprojlim \mathcal{H}_i$, along the projections $\{\theta_{j,i}\mid i<j\}$, is a valued hyperfield with value group $\Gamma$ and norm $\mathcal{N}:=\bigcup_i \mathcal{N}(\mathcal{H}_i)$. \label{Statement1}
        \item[(2)] \label{Statement2} $\mathcal{H}$ is stringent.
        \item[(3)] \label{Statement3}  For all $i\in I$, the factor hyperfield $\mathcal{H}_{\mathcal{N}(\mathcal{H}_{i})}$  is isomorphic to $\mathcal{H}_i$.
        \item[(4)] \label{Statement4} If the sequence $\bigl(\mathcal{N}_i(\mathcal{H}_i)\bigr)_{i\in I}$ tends to $\infty$, then $\mathcal{H}$ is a complete valued field.
    \end{enumerate}
\end{proposition}

\begin{proof}
As a set, the inverse limit of the given isometric system along its projections $\{\theta_{j,i}\mid i<j\}$ consists of all compatible sequences of elements of $\mathcal{H}_i$, $i\in I$, i.e.,
\[
\mathcal{H}=\left\{(a_i)_{i\in I}\in\prod_{i\in I}\mathcal{H}_i\mid \theta_{j,i}(a_j)=a_i\text{ for all }j>i\right\}.
\]

For  $a=(a_i)_i$, $b=(b_i)_i \in \mathcal{H}$, we set
    \begin{align}
        a+b:= \{ (c_i)_i \in \mathcal{H}\mid c_i \in a_i+b_i\text{ for all }i\in I\}
    \end{align}
    and
\begin{align}
    a b:= (a_i b_i)_i.
\end{align}
Moreover, we let $0:=(0_{\mathcal{H}_i})_i$ and  $1:=(1_{\mathcal{H}_i})_i$. Axiom (HF2) is trivially satisfied by $\mathcal{H}$. Regarding axiom (HF1), we note that all the axioms of canonical hypergroups are easy to verify, except associativity (CH1). Indeed, for $a=(a_i)_i$ and $b=(b_i)_i$ in $\mathcal{H}$, we have that
\begin{align*}
a+b& = \{ (c_i)_i \in \mathcal{H} \mid \text{ for } \ c_i \in a_i+b_i  \} \\
&= \{ (c_i)_i \in \mathcal{H} \mid \text{ for } \ c_i \in b_i+a_i  \}\\
 &= b+a.
\end{align*}
Set $-a:=(-a_i)_{i\in I}$. It is the unique additive inverse of $a$ as
\begin{align*}
0 \in a+b& \Longleftrightarrow \ 0_i\in a_i+b_i \in \mathcal{H} \ \text{ for all }i\in I; \\
& \Longleftrightarrow \  b_i=-a_i\ \text{ for all }i\in I \\
& \Longleftrightarrow \ b=-a.
 \end{align*}
 
\begin{align*}
a \in b+c& \Longleftrightarrow \ a_i\in b_i+c_i  \ \text{ for all }i\in I; \\
& \Longleftrightarrow \  c_i\in a_i-b_i \text{ for all }i\in I \\
& \Longleftrightarrow \ c\in a-b.
 \end{align*}

We showed (CH2), (CH3) and (CH4). Similarly, we show distributivity (HF3): 
\begin{align*}
(a+b)\cdot c & := \{ (l_i)_i\cdot (c_i)_i  \mid l_i \in a_i+b_i  \} \\
&= \{ (l_i')_i  \mid l_i' \in (a_i+b_i)\cdot c_i  \}\\
&= \{ (l_i')_i  \mid l_i' \in a_i c_i +b_i c_i  \}\\
 &=: a c + b c.
\end{align*}
     
It remains to show (CH1), that is, the associativity of $+$. For that, the valued structure is needed. We show first that $\mathcal{H}$ is naturally equipped with a map which behaves exactly as a valuation (although we cannot properly speak of valuations on $\mathcal{H}$ as we have not yet proved that $\mathcal{H}$ is a hyperfield). Recall that each $\mathcal{H}_i$ is equipped with an ultrametric $d_i$ as in Remark \ref{KrasnerValuedUltrametric}. We will now define an ultrametric $d$ on $\mathcal{H}$. Let $a=(a_i)_i$ and $b=(b_i)_i$. If for some $i\in I$, $a_i\neq b_i$, then for all $j>i$, we have that $a_j\neq b_j$. We want to show that in addition
\begin{align} \label{EquationDistanceEventuallyConstant}
d_j(a_j,b_j)= d_i(a_i,b_i).
\end{align}
Take $c_j\in a_j-b_j$, then $\val_j(c_j)=\val(\theta_{j,i}(c_j))$ 
by axiom (IH3). Since 
\[
\theta_{j,i}(c_j)\in \theta_{j,i}(a_j-b_j) \subseteq \theta_{j,i}(a_j)-\theta_{j,i}(b_j)=a_i-b_i
\]
and $a_i\neq b_i$, this shows that  $d_j(a_j,b_j)=\val_j(c_j)=d_i(a_i,b_i)$.

If $a\neq b$, then we define $d(a,b)$ to be the eventual value of $d_i(a_i,b_i)$. Otherwise, we set $d(a,b)$ to be $\infty$. We now consider the map
\begin{align*}
    \val:\mathcal{H}&\to\Gamma\cup\{\infty\}\\
    a&\mapsto d(a,0)
\end{align*}
Notice that, by what we have shown above, for all $a\in \mathcal{H}$, $\val(a)=\val_i(a_i)$ for all $i\in I$.\par 
We now show that $\val$ satisfies axioms (V0),(V1),(V2),(V3) and (V4). let $a=(a_i)_i$ and $b=(b_i)_i$ be arbitrary elements of $\mathcal{H}$.  By definition, $\val(a)=\infty$ if and only if $a=0=(0_{\mathcal{H}_i})_i$.
Moreover, for any $i\in I$, we have that
    \[\val(ab)=\val_i(a_ib_i)=\val_i(a_i)+\val_i(b_i)=\val(a)+\val(b).\]
    Thus, $\val$ satisfies axioms (V0) and (V1). If $c\in a+b$, then by axiom (V2) for $\mathcal{H}_i$, we have that $\val(c)=\val_i(c_i)\geq \min\{\val(a_i),\val(b_i)\}=\min\{\val(a),\val(b)\}$. Thus, $\val$ also satisfies (V2).
    If $0\notin a+b$, then clearly $0_{\mathcal{H}_i}\notin a_i+b_i$ for all $i\in I$. Now, fix some $i\in I$. If $c\in a+b$, then $\val(c)=\val_i(c_i)\in\val_i(a_i+b_i)$. Since by axiom (V3) applied to $\mathcal{H}_i$, we have that $\val(a_i+b_i)$ has a unique element, we may conclude that $\vert \val(a+b) \vert =1$ and so $\val$ satisfies (V3).
    Let us finally show (V4). Set $\mathcal{N}_i:=\mathcal{N}(\mathcal{H}_i)$ and $\mathcal{N}:= \bigcup_i \mathcal{N}_i$, thus $\mathcal{N}$ is the smallest initial segment which contains $\mathcal{N}_i$ for all $i\in I$. Take now $c,c'\in a+b$. Then $d(c,c')=d_i(c_i,c'_i)$ satisfies eventually 
    \[
    d_i(c_i,c_i') > \mathcal{N}_i+ \min\{\val(a_i),\val(b_i)\}\text{ for all }i\in I.
    \]
    In particular, $d(c,c')>\mathcal{N}+ \min\{\val(a),\val(b)\}$. Conversely, if $c\in a+b$ and $d(c,c')>\mathcal{N}+ \min\{\val(a),\val(b)\}$, then for all $i\in I$, we have that $d_i(c_i,c_i')=d(c,c') > \mathcal{N}_i+ \min\{\val(a_i),\val(b_i)\}$ and thus $c'_i\in a_i+b_i$ for all $i\in I$ so that $c'\in a+b$. We have now shown that $\val$ satisfies (V4) as well.

 We can now begin the proof of the associativity of the hyperoperation that we defined on $\mathcal{H}$. Recall that, by our assumption, $\mathcal{N}$ is the positive segment of a convex subgroup of $\Gamma$. We begin by proving the following claim:

\begin{claim}\label{ClaimUniqueness}
    Take $a,b,c \in \mathcal{H}$. There exists an element $l\in \mathcal{H}$ such that $l_i \in a_i+b_i+c_i$ for all $i\in I$. If, in addition, $0_{\mathcal{H}_i} \notin a_i+b_i+c_i$ for all $i\in I$, then this element is unique.
\end{claim}

\begin{proof}[Proof of Claim \ref{ClaimUniqueness}]

    Note first that applying (IH2'), we obtain that
    \[
    \theta_{j,i}(a_j+b_j+c_j) \subseteq \theta_{j,i}(a_j)+\theta_{j,i}(b_j)+\theta_{j,i}(c_j) = a_i+b_i+c_i
    \]
    for all $j>i$. We show the existential part of the claim. If for all $i$, $0_i \in a_i+b_i+c_i$,  we can set $l_i=0_{\mathcal{H}_i}$ for all $i\in I$. Assume that for some $i$, $0_i \in a_i+b_i+c_i$. By Proposition \ref{PropositionSumOfNElementsIsABall}, this means equivalently that 
    \begin{align}
        \val(x_i) \leq \mathcal{N}_i+\min(\val(a_i),\val(b_i),\val(c_i)). \label{EquationEquivalenceOintheSum}
    \end{align}
    for all $x_i\in a_i+b_i+c_i$.
     Since  for every $j>i$, $j\in I $,  $\val(a_j+b_j+c_j)=\val(a_i+b_i+c_i)$, a fortiori the same equality holds with $j>i$ and $x_j\in a_j+b_j+c_j$. Without lost of generality, we assume that $0_i \in a_i+b_i+c_i$ for every $i\in I$, and we fix some arbitrary $i\in I$ and let $n>i$ be such that $\mathcal{N}_n\geq 2 \mathcal{N}_i$. We will show that for all $x,y \in a_n+b_n+c_n $ we have that
    \begin{align}\label{EquationUniquenessInverseLimit} \theta_{n,i}(x)=\theta_{n,i}(y).    \end{align}
    Since, $\theta_{n,i}(x),\theta_{n,i}(y) \in a_i+b_i+c_i$, we have that
    \begin{align*}
        d_i(\theta_{n,i}(x),\theta_{n,i}(y))= d_n(x,y) & >  \mathcal{N}_n+ \min \{\val_n(a_n),\val_n(b_n),\val_n(c_n)\} \quad \text{by Proposition \ref{PropositionSumOfNElementsIsABall} $(ii)$}\\
        &\geq \mathcal{N}_i+(\mathcal{N}_i+\min \{\val_i(a_i),\val_i(b_i),\val_i(c_i)\}) \\
        & \geq \mathcal{N}_i+\min\{ \val_i(\theta_{n,i}(x)),\val_i(\theta_{n,i}(y))\} \quad \text{by (\ref{EquationEquivalenceOintheSum})}.
    \end{align*}

  Using axiom (V4), we conclude that $0_{\mathcal{H}_i} \in \theta_{n,i}(x)-\theta_{n,i}(y)$ and thus \eqref{EquationUniquenessInverseLimit} must hold.\par 
  We define $l_i:=\theta_{n,i}(l_n)$ where $n>i$ is such that $\mathcal{N}_n\geq 2 \mathcal{N}_i$ and $l_n\in a_n+b_n+c_n$. This is clearly independent from the choice of $n$ and we have just shown that it is also independent from the choice of $l_n$.\par 
  This procedure has been described for an arbitrary $i\in I$ and yields a sequence $(l_i)_{i\in I}$ which satisfies $l_i\in a_i+b_i+c_i$ for all $i\in I$. By construction, it follows that $\theta_{j,i}(l_j)=\theta_{j,i}(\theta_{n,j}(l_n))=\theta_{n,i}(l_n)=l_i$
  for all $i<j$, where $n>j>i$ has been chosen such that $\mathcal{N}_n>2\mathcal{N}_j\geq\mathcal{N}_i$. Hence, $(l_i)_{i\in I}\in\mathcal{H}$ and we have proved the existence part of our claim. The uniqueness part follows immediately from \eqref{EquationUniquenessInverseLimit}.\qedhere
\end{proof}
    
To prove associativity, take three elements $a,b,c\in \mathcal{H}$. We will prove that $(a+b)+c\subseteq a+(b+c)$. The other inclusion can be treated symmetrically. Pick $l\in (a+b)+c$. We need to show that $l\in a+(b+c)$. We distinguish two cases:\par
\textbf{Case 1:} Assume that $0_{\mathcal{H}_i} \in a_i+b_i+c_i$ for all $i\in I$. This means that $-a_i\in b_i+c_i$ and $-c_i\in a_i+b_i$ for all $i\in I$, thus by the definition of the hyperoperation in $\mathcal{H}$, we have that $-a\in b+c$ as well as $-c\in a+b$.

\textbf{Case 1.1:} If $\val(a)>\mathcal{N} +\min\{\val(b),\val(c)\}$, then

\begin{center}
\begin{tikzpicture}
\draw [decorate,decoration={brace,amplitude=5pt,raise=1ex}](3.5,0) -- (6.5,0) node[midway,yshift=1.5em]{$\mathcal{N}$};
\draw (3.5,0) -- (6.5,0);
\draw (3.5,0) node[left]{$a$} ;
\draw (0,-1) -- (3,-1); 
\draw [decorate,decoration={brace,amplitude=5pt,raise=1ex}](0,-1) -- (3,-1) node[midway,yshift=1.5em]{$\mathcal{N}$};
\draw (0,-1) node[left]{$b$} ;
\draw (0,-2) -- (3,-2);
\draw (0,-2) node[left]{$c$} ;
\end{tikzpicture}
\end{center}

we claim that $c=-b$ in this case. Indeed, since $\val_i(0+a_i)>\mathcal{N} +\min\{\val_i(b_i),\val_i(c_i)\}$ and $-a_i\in b_i+c_i$, by the norm axiom $0\in b_i+c_i$ and $c_i=-b_i$ for all $i$.
We show that $(a+b)+c= b-b= a+(b+c)$. Let $e\in a+b$ such that $l\in e+c$. By Remark \ref{RemarkSumIsSingleton}, for all $i$, $a_i+b_i=\{b_i\}$. It follows that $e=b$ and $l \in b+c = b-b$. By Axiom (V2), $\val(l_i-a_i-0)\geq \min(\val(l_i), \val(a_i)) > \mathcal{N} +\val(b_i) $ and thus by (V4) $l_i-a_i\subseteq b_i-b_i$. Let $x\in l-a$. By Axiom $\text{(CH4)}_\mathcal{H}$ that we proved earlier, we have that $l \in a + x$ and $l\in a+(b+c)$.  

\textbf{Case 1.2:} If $\val(a) \leq\mathcal{N} +\min\{\val(b),\val(c)\}$, then
\begin{center}
\begin{tikzpicture}
\draw [decorate,decoration={brace,amplitude=5pt,raise=1ex}](1.5,0) -- (4.5,0) node[midway,yshift=1.5em]{$\mathcal{N}$};
\draw (1.5,0) -- (4.5,0);
\draw (1.5,0) node[left]{$a$} ;
\draw (0,-1) -- (4.5,-1); 
\draw [decorate,decoration={brace,amplitude=5pt,raise=1ex}](0,-1) -- (4.5,-1) node[midway,yshift=1.5em]{$\mathcal{N}$};
\draw (0,-1) node[left]{$b$} ;
\draw (0,-2) -- (4.5,-2);
\draw (0,-2) node[left]{$c$} ;
\end{tikzpicture}
\end{center}

since by assumption $(\mathcal{N}_i)_{i\in I}$ is cofinal in a convex subgroup of $\Gamma$, we have that 
\begin{equation}\label{cutrightequal}
(\Gamma_{<0}\cup\mathcal{N})+\val(a)=(\Gamma_{<0}\cup\mathcal{N})+\min\{\val(b),\val(c)\}.
\end{equation}
Fix $i\in I$. By Proposition \ref{PropositionSumOfNElementsIsABall} $(ii)$ and (V2), since $-a_i\in b_i+c_i$, we have that 
\[
\val(l)=\val_i(l_i)>\mathcal{N}_i+\min\{\val_i(a_i),\val_i(b_i),\val_i(c_i)\}=\mathcal{N}_i+\min\{\val_i(b_i),\val_i(c_i)\}.
\]
By the fact that $i$ was arbitrary, we conclude that $\val(l)>\mathcal{N}+\min\{\val(b),\val(c)\}$ and follows from \eqref{cutrightequal} that $\val(l)>\mathcal{N}+\val(a)$. Hence, from (V4) we obtain that $l\in a-a\subseteq a+(b+c)$.\par
\textbf{Case 2:} If $0_{\mathcal{H}_i}\notin a_i+b_i+c_i$ for some $i\in I$, then as we have done above we can assume that this happens for all $i\in I$. Pick any element $e\in a+(b+c)$ (which exists by Claim \ref{ClaimUniqueness}). By definition we have that $e_i\in a_i+b_i+c_i$ for all $i\in I$. On the other hand, $l_i\in a_i+b_i+c_i$ for all $i\in I$ as well. We conclude by the uniqueness part of Claim \ref{ClaimUniqueness} that $l=e\in a+(b+c)$ as we wished to show. 

We have now finally shown that $\mathcal{H}$ is a valued hyperfield. By Claim \ref{ClaimUniqueness} it immediately follows that $\mathcal{H}$ is stringent.
    
We now begin the proof of part \textit{(3)}. Fix $i\in I$ and observe that the map $\theta_i: \mathcal{H} \rightarrow \mathcal{H}_i$,  $a\mapsto a_i$ is an isometric homomorphism, indeed that it satisfies (IH1) and (IH3) can be shown easily: $a,b\in \mathcal{H}$ $\theta_i(a b)=a_ib_i=\theta_i(a) \theta_i(b)$, and $\val(a)=\val_i(\theta_i(a))$ by definition.
    We now aim to show that (IH2) holds. Since $\theta_i(0)=0_i$ and $\mathcal{H}_i$ is a valued hyperfield, by Lemma \ref{PropositionIH2EquivalentIH2'} it suffices to show that (IH2') holds for $\theta_i$. This follows immediatly from the definitions of $\theta_i$ and the hyperoperation of $\mathcal{H}$: if $l=(l_i)\in a+b$, then $\theta_i(l)=l_i \in a_i+b_i=\theta_i(a)+\theta_i(b)$.
    By Lemma \ref{LemmaIsometryIndicesIsomophism}, the map
    \begin{align*}
    \phi:\mathcal{H}_{\mathcal{N}_i}&\to\mathcal{H}_i\\
    [a]_{\mathcal{N}_i}&\mapsto a_i
    \end{align*}
    is an isomorphism of valued hyperfields.
    
    We now come to the proof of part \textit{(4)}. We then assume that $\mathcal{N}=\Gamma_{\geq0}$.By the norm axiom (V4), we have that if $a,b\in \mathcal{H}$, then for all $c,c'\in a+b$ we have that
    \[
    \val(x) > \mathcal{N}+\min\{\val(a),\val(b)\}
    \] 
    for all $x\in c-c'$. This can only be true if $c-c'=\{0\}$, so in particular $c=c'$ must hold. This shows that $\mathcal{H}$ is a valued field in this case.\par
    Let us show that it is complete. Let $(a^\nu)_{\nu<\lambda}$ be a Cauchy sequence in $\mathcal{H}$: for all $\gamma\in\Gamma$ there exists $N<\lambda$ such that $\val(a^\nu-a^\mu)>\gamma$ for all $\nu,\mu\geq N$.
    If $\bigl(\val(a^\nu)\bigr)_{\nu<\lambda}$ is not bounded from above in $\Gamma$, then $(a^\nu)_{\nu<\lambda}$ tends to $0$. Otherwise, let $\delta\in\Gamma$ be an upper-bound for $(\val(a^\nu))_{\nu<\lambda}$. For all $i\in I$, let $\nu_i<\lambda$ be such that $\val(a^\nu-a^\mu)> \mathcal{N}_i+\delta$ for all $\nu,\mu\geq\nu_i$. 
    For $i\in I$, set $a_i:=a_{i}^{\nu_i}$. We now show that $(a_i)_{i\in I}$ is an element of $\mathcal{H}$. 
    For all $j>i$, we have that $\theta_{j,i}(a_j)=\theta_{j,i}(a^{\nu_j}_j)=a^{\nu_j}_i$. If $a^{\nu_j}_i\neq a^{\nu_i}_i=a_i$, then
    \[
    d_i(a^{\nu_j}_i,a^{\nu_i}_i)=d(a^{\nu_j},a^{\nu_i})>\mathcal{N}_i+\delta\geq\mathcal{N}_i+\min\{\val(a^{\nu_i}),\val(a^{\nu_j})\}=\mathcal{N}_i+\min\{\val_i(a^{\nu_i}_i),\val_i(a^{\nu_j}_i)\}.
    \]
    It follows that $0_{\mathcal{H}_i}\in a^{\nu_j}_i-a^{\nu_i}_i$. This contradiction shows that $\theta_{j,i}(a_j)=a_i$ must hold. Thus, $a:=(a_i)_{i\in I}$ is an element of $\mathcal{H}$.

    We can now easily see that it is the limit of the Cauchy sequence $(a^\nu)_{\nu<\lambda}$. Indeed, since $\mathcal{N}=\Gamma_{\geq0}$, by definition we have that $\bigl(\val(a^\nu-a)\bigr)_{\nu<\lambda}$ tends to $\infty$. This completes the proof of the proposition.  
\end{proof}

Let us now describe an interesting example.

\begin{example}\label{ExamplesProjectiveLimiteHyperfields} 
Consider the field $\mathcal{K} := \mathbb{Q}(x)(y)$ equipped with the composition $\val:= \val_x \circ \val_y : \mathcal{K} \rightarrow \mathbb{Z}\times\mathbb{Z}$ of the is the $x$-adic valuation $\val_x$ and the $y$-adic valuation $\val_y$. Explicitly, we have for $a:= \sum_n Q_n(x)y^n\in \mathcal{K}$ where $Q_n(x)\in \mathbb{Q}(x)$ that  
\[
\val\left(\sum_n Q_n(x)y^n\right):= (\val_y(a), \val_x(Q_{\val_y(x)})).
\]
Then an element in $\mathcal{H}_{(0,n)}$ are all of the form $y^mx^k\sum_{i\leq n} a_ix^i + \mathfrak{m}_{(0,n)}$, where $a_i\in \mathbb{Q}$ for $i\leq n$ (consider the euclidean division according the increasing power). 
    
The sequence $(\mathcal{H}_{(0,n)})_n$ is a isometric system of hyperfields and by Proposition \ref{PropositionInverseLimitHyperfields}, the inverse limit $\varprojlim_{n}\mathcal{H}_{(0,n)}$ is an stringent hyperfield. We see that it is isomorphic to the valued hyperfield $\mathcal{H}_{\Delta}(\mathcal{K}')$ where $\mathcal{K}'=\mathbb{Q}((x))((y))$ is the extension of $\mathcal{K}$ to the completion for the $x$-adic valuation  and $\Delta := \{0\}\times \mathbb{N}$.

    Indeed, an element of the inverse limit $\varprojlim_{n}\mathcal{H}_{(0,n)}$ is of the form
    \[\left(y^mx^k\sum_{i<n} a_ix^i \right)_{n\in \mathbb{N}_{>0}},\]
    where $m,k \in \mathbb{Z}$ and $(a_i)_{i\in \mathbb{N}} \in \mathbb{Q}^\mathbb{N}$, and an element of $\mathcal{H}_{\Delta}(\mathcal{K}')$ is of the form
    \[y^m \sum_{i}a_ix^{k+i}\]
    where again $m,k \in \mathbb{Z}$ and $(a_i)_{i\in \mathbb{N}} \in \mathbb{Q}^\mathbb{N}$; the map 
        \[\left(y^mx^k\sum_{i<n} a_i       x^i \right)_{n\in \mathbb{N}_{>0}} \longmapsto y^m \sum_{i}a_ix^{k+i},\]
        gives an isometric isomorphism of hyperfields between
        $\varprojlim_{n}\mathcal{H}_{(0,n)}$ and $\mathcal{H}_{\Delta}(\mathcal{K}')$.
        
Similarly, by Proposition \ref{PropositionInverseLimitHyperfields}, the inverse limit $\varprojlim_{n}\mathcal{H}_{(n,0)}$ is a valued field, and it is isomomorphic to  the $y$-adic completion $\mathbb{Q}(x)((y))$ of $\mathcal{K}$.

Note that the inverse limit
$\varprojlim_{m}\mathcal{H}_{(1,m)}$ is not a hyperfield when equipped with the addition $+$ as defined in the proof of Proposition \ref{PropositionInverseLimitHyperfields}. 
    
    For example, the set 
    \[(\sum_{i<n} x^i)_n + (y\sum_{i<n} x^i)_n := \left\{(c_n)_n\in \varprojlim_{n} \mathcal{H}_{(1,n)} \ \vert \ c_n\in \sum_{i<n} x^i +(y\sum_{i<n} x^i) \right\}\]
    is empty, as no compatible sequence $(c_n)_n$ can satisfy $c_n \in \sum_{i<n} x^i +(y\sum_{i<n} x^i) $ for every $n$.
\end{example}

    \section{On stringent valued hyperfields}\label{SectionAxiomatisationStringentHyperfields}
    
We give a useful axiomatisation of stringent valued hyperfields, which describes these structure with only standard (single-valued) operations. This will be useful in the next section when we construct a henselian valued field from a stringent valued hyperfield. Part of the results of this section can be found in \cite{Touthesis}.

        \begin{definition}
        A \emph{sequence-structure} is a tuple $(\mathcal{F}^\times,\mathcal{H}^\times,\Gamma,\iota,\nu)$ where $\mathcal{F}$ is a hyperfield, $\mathcal{H}^\times$ an abelian group and $\Gamma$ an ordered abelian group and
        \[
        \begin{tikzcd}
        \{1\}\arrow[r]& \mathcal{F}^\times \arrow[r,"\iota"]& \mathcal{H}^\times\arrow[r,"\nu"]& \Gamma\arrow[r]& \{0\}.
        \end{tikzcd}
        \]
        is a short-exact sequence of abelian groups.
        \end{definition}
        Given a sequence-structure $(\mathcal{F}^\times,\mathcal{H}^\times,\Gamma,\iota,\nu)$ we will always identify $\mathcal{F}^\times$ with its image under~$\iota$~in~$\mathcal{H}^\times$.

\begin{theorem}[{\cite[Theorem 1.2 ]{BS21}}]\label{BowlerSu}
    A hyperfield $\mathcal{H}$ is stringent if and only if it arises from a sequence structure $(\mathcal{F}^\times,\mathcal{H}^\times,\Gamma,\iota,\nu)$ where $\mathcal{F}$ is either $\mathbb{K}$, $\mathbb{S}$ or a field.
\end{theorem}
        
        Let us recall some details of the construction described by Bowler and Su in \cite{BS21}.\par

        Given a sequence-structure $(\mathcal{F}^\times,\mathcal{H}^\times,\Gamma,\iota,\nu)$, complete the embedding of $\mathcal{F}^\times$ in $\mathcal{H}^\times$ by adding a new element $\mathbf{0}$ to $\mathcal{H}^\times$, absorbing for the multiplication (\textit{i.e.} $\mathbf{0}\cdot \mathbf{a} =\mathbf{0}$ for all $\mathbf{a}\in\mathcal{H}^\times \cup \lbrace \mathbf{0}\rbrace$). Denote by $\mathcal{H}$ the set $\mathcal{H}^{\times} \cup \lbrace \mathbf{0}\rbrace$ and then add a new element $\infty$ to $\Gamma$, with $\infty>\Gamma$. Extend the map $\nu$ to $\mathcal{H}$ by setting $\nu(\mathbf{0}):=\infty$. Now, equip $\mathcal{H}$ with the following hyperoperation $\boxplus$:
       \[
            \mathbf{a} \boxplus \mathbf{b} = \begin{cases}
            \iota(1-1) \cdot \mathbf{a}  \cup \{  \mathbf{x}\mid \nu(\mathbf{x})> \nu(\mathbf{a})  \}  & \text{ if } \mathbf{b}=(-1) \cdot \mathbf{a},\\
            \{\mathbf{a}\} & \text{ if } \nu(\mathbf{b})>\nu(\mathbf{a}),\\
            \{\mathbf{b}\} & \text{ if } \nu(\mathbf{a})>\nu(\mathbf{b}),\\
            \{(\mathbf{a}/\mathbf{b}+1)\cdot \mathbf{b}\}&  \text{ otherwise.}
        \end{cases}\quad\quad\quad(\mathbf{a},\mathbf{b}\in\mathcal{H})
        \]
Notice that to define $\boxplus$ we used addition and subtraction only in $\mathcal{F}=\ker\nu$. In \cite[Lemma 3.1 \& Lemma 4.2]{BS21} it is shown that with this hyperoperation $\mathcal{H}$ becomes a hyperfield. We denote by $\mathbf{a}\boxminus\mathbf{b}$ the sum $\mathbf{a}\boxplus(-\mathbf{b})$, where $\mathbf{a}, \mathbf{b}\in \mathcal{H}$. It is clear from the definition of $\boxplus$ that, if $\mathcal{F}$ is either $\mathbb{K}$ or $\mathbb{S}$ or a field, then the resulting hyperfield $\mathcal{H}$ is stringent.

Conversely, given a stringent hyperfield $(\mathcal{H},\boxplus,\cdot,\mathbf{0},\mathbf{1})$ an ordered abelian group $\Gamma$ can be obtained as the quotient of $\mathcal{H}^\times$ modulo the equivalence relation
\[
a\sim b\Longleftrightarrow a\boxplus b\neq \{a\} \wedge a\boxplus b \neq \{b\}.
\]
In \cite[Lemma 3.10]{BS21}, Bowler and Su show that $\mathcal{F}:=\ker\nu\cup\{\mathbf{0}\}$, where $\nu:\mathcal{H}^\times\to\Gamma$ is the canonical epimorphism, is either $\mathbb{K}$ or $\mathbb{S}$ or a field. 

It follows from Theorem \ref{BowlerSu} that a stringent valued hyperfield $(\mathcal{H},\val)$ always arise from a sequence-structure $(\mathcal{F}^\times,\mathcal{H}^\times,\Gamma,\iota,\nu)$, where $\mathcal{F}$ is a field. Indeed, in $\mathbb{K}$ and $\mathbb{S}$ we have $1\in 1 - 1$ and thus if $\mathcal{F}\in\{\mathbb{K},\mathbb{S}\}$, then $\textbf{1}\in \textbf{1}\boxminus\textbf{1}$ would be true in $\mathcal{H}$. On the other hand, this would prevent $\mathcal{H}$ from admitting a valuation map $\val$ in the sense of Definition \ref{DefinitionValuationHyperfields} as then $\val(\textbf{1})=0$ while $\val(\textbf{x})>0$ for all $\textbf{x}\in \textbf{1}-\textbf{1}$ must hold by the norm axiom.
Conversely, it is not difficult to show that if $(\mathcal{F}^\times,\mathcal{H}^\times,\Gamma,\iota,\nu)$ is a sequence-structure and if $\mathcal{F}$ is a field, then the map $\nu$ is a valuation, making $\mathcal{H}$ a valued hyperfield of norm \{0\}.

Let us now observe that stringent valued hyperfields can be identified (bi-interpretable in the sense of model theory) with a structure $(\mathcal{H},\oplus,\cdot,0,1)$ where $\oplus$ is a binary operation: there is no need for the concept of hyperoperation in this case.  

\begin{definition}\label{DefinitionRVSorts}
               Let $\mathcal{H}$ a stringent valued hyperfield given by a sequence-structure $(\mathcal{F}^\times,\mathcal{H}^\times,\Gamma,\iota,\val)$
        For two elements $\mathbf{a}$ and $\mathbf{b}$ in $\mathcal{H}$, we denote by $\textbf{a} \oplus \textbf{b}$ the element:
        \[\textbf{a} \oplus \textbf{b} = \begin{cases}
            \textbf{0} & \text{ if } \textbf{\textbf{a}}=\textbf{\textbf{b}}=0,\\
            \textbf{a} & \text{ if } \val(\textbf{b})>\val(\textbf{a}),\\
            \textbf{b} & \text{ if } \val(\textbf{a})>\val(\textbf{b}),\\
            (\textbf{a}/\textbf{b}+1)\cdot \textbf{b} & \text{ otherwise.}
        \end{cases}
        \]
\end{definition}

We give below a list of axioms for the structure $(\mathcal{H},\oplus,\cdot,\mathbf{0},\mathbf{1})$ which arises from a stringent valued hyperfield. Notice that $\oplus$ extends by definition the addition in $\mathcal{F}$ and that unlike the hyperaddition $\boxplus$, it is not (always) associative. As a consequence, we have to adopt a convention for $\bigoplus\limits_{i<n} \mathbf{a}_i$:
    \begin{remark}\label{RemarkConventionSumOplus}
        Let $n\in \mathbb{N}$ and $\textbf{a}_1, \ldots, \mathbf{a}_n \in \mathcal{H}$. We will use the notation $\bigoplus\limits_{i<n} \textbf{a}_i$ for $((((\textbf{a}_{1} \oplus \textbf{a}_{2})\oplus \textbf{a}_{3})\oplus \cdots )\oplus \textbf{a}_{n})$ only when the sum $(((\textbf{a}_{\sigma(1)} \oplus \textbf{a}_{\sigma(2)})\oplus \textbf{a}_{\sigma(3)}) + \cdots )\oplus \textbf{a}_{\sigma(n)})$ is associative for any permutation $\sigma$ (in other words, when this term does not depend on the position of the parenthesis). It happens in particular when all $\textbf{a}_i$'s have the same valuation (dividing by $\textbf{a}_1$, we get a sum in the field $\mathcal{F}$) or \textit{a contrario}, when they have pairwise distinct valuations (by additive absorption). 
      \end{remark}

    We describe the essential properties of $\oplus$. Let $\mathrm{L}$ be the one-sorted language with signature $\lbrace\oplus,\cdot,\textbf{0},\textbf{1}\rbrace$. We will see that hyperfield equipped with this operation $\oplus$ are exactly the $\mathrm{L}$-structures satisfying the following list of axioms (RV1-RV7): 
        \begin{enumerate}
            \item[(RV1)] $(\mathcal{H}^\times,\cdot,\textbf{1})$ is an abelian group, where $\mathcal{H}^\times=\mathcal{H} \setminus \lbrace \textbf{0}\rbrace$,
            \item[(RV2)] (neutral element for $\oplus$) $\forall \textbf{a}\in \mathcal{H}, \ \textbf{0} \oplus \textbf{a} = \textbf{a}, $
            \item[(RV3)] (semi or half-associativity) $[(\textbf{a}\oplus \textbf{b}) \oplus \textbf{c} \neq \textbf{a} \oplus (\textbf{b} \oplus \textbf{c})] \Rightarrow \textbf{a}\oplus \textbf{b} =\textbf{0} \text{ or } \textbf{b} \oplus \textbf{c}=\textbf{0},$
            \item[(RV4)] (commutativity for $\oplus$) $\forall \textbf{a},\textbf{b} \in \mathcal{H}, \ \textbf{a}\oplus \textbf{b} =\textbf{b} \oplus \textbf{a}$,
            \item[(RV5)] (distributivity) $\forall \textbf{a},\textbf{b},\textbf{c} \in \mathcal{H}, \ (\textbf{a} \oplus \textbf{b}) \cdot \textbf{c} = \textbf{ac} \oplus \textbf{bc} $.
        \end{enumerate}
        We define $\mathcal{F}^\times$ as the set:
        $$\lbrace \textbf{r}\in \mathcal{H} \setminus\lbrace \textbf{0}\rbrace \ \vert \ \textbf{1} \oplus \textbf{r} \neq \textbf{1} \ \text{ and } \ \textbf{1}\oplus \textbf{r}^{-1} \neq \textbf{1} \rbrace.$$
        We write $\mathcal{F}:=\mathcal{F}^\times \cup \lbrace 0 \rbrace$ and we may denote its elements $a,b,c,\ldots\in \mathcal{F}$ with the usual font and denote the restriction $\restriction{\oplus}{\mathcal{F}}$ by the symbol $+$.
        \begin{enumerate}
            \setcounter{enumi}{5}
            \item[(RV6)] $(\mathcal{F}:=\mathcal{F}^\times \cup \lbrace 0 \rbrace,\cdot,+,0,1)$ is a field, 
            \item[(RV7)] (uniform additive absorption) $\forall \textbf{a}\in \mathcal{H}, \ \forall r \in \mathcal{F}^\times, \ (\textbf{a} \oplus \textbf{1} = \textbf{1}) \Longleftrightarrow (\textbf{a} \oplus r = r)$.
        \end{enumerate}
        
        We see easily that any hyperfield $\mathcal{H}$ equipped with the low $\oplus$ as defined above satisfies the axioms (RV1-RV7). We derive from them the property (RV8-RV10). This shows that any structure $(\mathcal{H},\oplus,\cdot,\textbf{0},\textbf{1})$ satisfying (RV1-RV7) rise from a sequence-structure $1\rightarrow \mathcal{F}^\times \rightarrow \mathcal{H}^\times \rightarrow \Gamma \rightarrow 0$ where $\mathcal{F}$ is a field. By the previous paragraph, this means that $\mathcal{H}$ is a stringent valued hyperfield.
        \begin{enumerate}
            \setcounter{enumi}{7}
            \item[(RV8)] (multiplicative absorption) $\forall \textbf{a}\in \mathcal{H}, \ \textbf{0} \cdot \textbf{a} = \textbf{0}. $
        \end{enumerate}
            Assume for some $\textbf{a}\in \mathcal{H}$, $\textbf{0}\cdot \textbf{a} \neq \textbf{0}$. Then as $\textbf{0}$ is a neutral element for $\oplus$ and by distributivity $\textbf{0}\cdot \textbf{a} = \textbf{0}\cdot \textbf{a} \oplus \textbf{0}\cdot \textbf{a}$. We may multiply by the multiplicative inverse of $\textbf{0}\cdot \textbf{a}$ and we get a contradiction in $\mathcal{F}$.
        \begin{enumerate}
            \setcounter{enumi}{8}
            \item[(RV9)] (additive inverse)  $\forall \textbf{a} \in \mathcal{H}, \exists ! \textbf{b} \ \textbf{a}\oplus \textbf{b} = \textbf{0}$.   
        \end{enumerate}
        This inverse is given by $-\textbf{a} := -\textbf{1}\cdot \textbf{a}$. Indeed we have $(\textbf{a} \oplus -\textbf{a})=\textbf{a}\cdot(\textbf{1}+-\textbf{1})=\textbf{0}$. Uniqueness is clear if $\textbf{a}=\textbf{0}$. Assume $\textbf{a}\neq \textbf{0}$, if $\textbf{b}\in \mathcal{H}^\times$ is such that $\textbf{a}\oplus \textbf{b}=\textbf{0}$, in particular $\textbf{b} \neq \textbf{0}$. Then by distributivity and multiplicative absorption $\textbf{b}/\textbf{a} \in k^\times$ and from $\textbf{b}/\textbf{a}+\textbf{1}=\textbf{0}$ we get $\textbf{b}=-\textbf{a}$.\\
        
        We recover the value group by setting $\Gamma := \mathcal{H}^\times/\mathcal{F}^\times$. For the order in $\Gamma$, one must define it as follows:
        \[
        \forall [\textbf{a}],[\textbf{b}] \in \Gamma, [\textbf{a}] < [\textbf{b}] \ \Longleftrightarrow \ \textbf{1}\oplus \textbf{b}/\textbf{a}= \textbf{1} \ \Longleftrightarrow \ \textbf{a}\oplus \textbf{b}= \textbf{a},
        \]
        where $[\textbf{a}]$ denote the class of $\textbf{a}$ modulo $\mathcal{F}^\times$. By uniform additive absorption and distributivity, this definition does not depend of the representative $\textbf{a}$ and $\textbf{b}$ we have chosen. Indeed, if $r,r' \in \mathcal{F}^\times$, one gets:
        
        \[
        1 \oplus \textbf{b}r/\textbf{a}r'=\textbf{1} \Longleftrightarrow r'/r \oplus \textbf{b}/\textbf{a}= r'/r \overset{(7)}{\Longleftrightarrow} \textbf{1}\oplus \textbf{b}/\textbf{a} =\textbf{1}.
        \]
        \begin{enumerate}
            \setcounter{enumi}{9}
            \item[(RV10)] $(\Gamma,<)$ is an ordered group.   
        \end{enumerate}
        
        Anti-symmetry of $<$ follows from the definitions of $\mathcal{F}^\times$ and $<$, and transitivity is given by semi-associativity:
        Assume $\textbf{a}\oplus \textbf{b}=\textbf{a}$ and $\textbf{b} \oplus \textbf{c}=\textbf{b}$, then either $\textbf{b}=\textbf{0}$, or $\textbf{a} \oplus \textbf{b} \neq \textbf{0}$ and $\textbf{b} \oplus \textbf{c} \neq \textbf{0}$. In any case, $\textbf{a} \oplus \textbf{c} = (\textbf{a} \oplus \textbf{b}) \oplus \textbf{c} = \textbf{a} \oplus (\textbf{b} \oplus \textbf{c}) = \textbf{a} \oplus \textbf{b} = \textbf{a}$. It's a total order since for all $\textbf{a},\textbf{b} \in \mathcal{H}^\times$, either $\textbf{a}/\textbf{b} \in \mathcal{F}^\times$, $\textbf{a}/\textbf{b}\oplus \textbf{1} = \textbf{1}$ or $\textbf{b}/\textbf{a} \oplus \textbf{1} =\textbf{1}$, which respectively gives $[\textbf{a}]=[\textbf{b}]$, $[\textbf{a}]>[\textbf{b}]$ or $[\textbf{b}]>[\textbf{a}]$. 
        We complete the valuation map by setting $\val(\textbf{0}) =\infty$ where $\infty> \Gamma$            \footnote{To avoid the use of conventions for $\textbf{0}$, one might define the quotient $\mathcal{H}/\mathcal{F}^\times$ as the set of orbits of $\mathcal{H}$ under the action of $\mathcal{F}^\times$. This action preserve the multiplication in $\mathcal{H}$. We get that $\textbf{0}$ is the unique element in its orbit $[\textbf{0}]$, which we denoted by $\infty$. Then the definition of the ordering $<$ gives that $\infty> \Gamma$.}.

  \begin{definition} 
      
        We call \emph{$\RV$-sort} any structure $\mathcal{H}$ in the language $\{\oplus, \cdot, \mathbf{0},\mathbf{1}\}$ satisfying the axioms (RV1-RV7). 
  
  \end{definition}
  
    We sum up the three equivalent points of view for stringent valued hyperfields in the following fact:  
    \begin{fact}\label{biinterpretability}
        There is a one-to-one correspondence between sequence-structures, stringent hyperfields and $\RV$-sorts.  The following structures are pairwise bi-interpretable:
        \begin{itemize}
            \item[] (Sequence-structure) $\{(\mathcal{H}^\times,\cdot,1), (\mathcal{F}^\times,+,\cdot,0,1), (\Gamma,0,+,<), \iota, \val \}$, where $\mathcal{F}$ is a field.
            \item[] (Stringent valued hyperfield) $(\mathcal{H}, \boxplus, \cdot, 0, 1)$ (where the hyperoperation $\boxplus$ is encoded by the ternary predicate $z\in x\boxplus y$).
            \item[] ($\RV$-sort) $(\mathcal{H}, \oplus, \cdot, 0, 1)$ (where $\oplus$ is a binary function).
        \end{itemize}
    \end{fact}

\section{The Hahn field of a stringent valued hyperfield}\label{SectionHahnProductStingentHyperfield}

    Let $\mathcal{H}$ be a stringent valued hyperfield, with value group $\Gamma$ and valuation $\val_\mathcal{H}$. We give a construction, based on the so called Hahn-product for valued fields, which produces a henselian valued field $\mathcal{K}$ such that $\mathcal{H}_{\{0\}}(\mathcal{K})\simeq\mathcal{H}$. By Fact \ref{biinterpretability} we can (and we will) consider the hyperfield $\mathcal{H}$ as an RV-sort, that is, equipped with the structure $(\mathcal{H}, \oplus,\cdot, \mathbf{0},\mathbf{1})$. This point of view seems to be more adequate in order to treat this problem.

    \begin{definition}[The Hahn field $\mathcal{H}^{(\Gamma)}$] \index{Hahn field associated to an $\mathcal$-sort} 
        The \emph{Hahn field} $\mathcal{H}^{(\Gamma)}$ associated to the $\RV$-sort $(\mathcal{H},\oplus,\cdot, \mathbf{0},\mathbf{1})$ is the set: 
        \[ \left\lbrace (\textbf{a}_\gamma)_{\gamma \in \Gamma} \ \vert \ \forall \gamma\in \Gamma \ \mathbf{a}_\gamma \in \mathcal{H}, \val_{\mathcal{H}}(\mathbf{a}_\gamma) \in \lbrace \gamma,\infty \right\rbrace
         \text{ and } \supp(\mathbf{a}_\gamma)_{\gamma} \ \text{is well-ordered} \rbrace\]
        where $\supp(\mathbf{a}_\gamma)_\gamma = \lbrace \gamma \in \Gamma \ \vert \ \mathbf{a}_{\gamma} \neq \textbf{0}\rbrace$.\par
        We equip $\mathcal{H}^{(\Gamma)}$ with the following law:
        \[
        (\mathbf{a}_\gamma)_\gamma + (\mathbf{b}_\gamma)_\gamma := (\mathbf{a}_\gamma \oplus \mathbf{b}_\gamma)_\gamma,
        \]
        \[
        (\mathbf{a}_\gamma)_\gamma \cdot (\mathbf{b}_\gamma)_\gamma := \left(\bigoplus\limits_{\delta+\epsilon =\gamma}\mathbf{a}_{\delta}\cdot \mathbf{b}_{\epsilon}\right)_\gamma,
        \]
        and with the constant $0:=(\mathbf{0})_\gamma$ and $1:= (\mathbf{a}_\gamma)_\gamma\in \mathcal{H}^{(\Gamma)}$  where 
        \[\mathbf{a}_\gamma= 
            \begin{cases} 
                \textbf{0} \text{ if } \gamma \neq 0 \\
                \textbf{1} \text{ if } \gamma=0.
            \end{cases}\]

        Notice that the product is well-defined as the supports of $a$ and $b$ are well-ordered, the term $\bigoplus\limits_{\delta+\epsilon =\gamma}\mathbf{a}_{\delta}\cdot \mathbf{b}_{\epsilon}$ (for a fix $\gamma$) is finite. Moreover, since elements have valuation $\gamma$ or $\infty$, the sum is independent of the  choice of parenthesis (see Remark \ref{RemarkConventionSumOplus}).  Then $\supp(a+b) \subset \supp(a)+\supp(b)$ and in particular, $\supp(a+b)$ is a well-ordered set of $\Gamma$ and $a\cdot b \in \mathcal{H}^{(\Gamma)}$. We set 
        \[
        \val(\mathbf{a}_\gamma)_\gamma := \min \supp(\mathbf{a}_\gamma)_\gamma.
        \]
    \end{definition}
    with the convention that $\val(0)=\infty$.  An element $(\mathbf{a}_\gamma)_{\gamma \in \Gamma} \in \mathcal{H}^{(\Gamma)}$ is written $\sum\limits_{\gamma \in \Gamma} \mathbf{a}_\gamma$ where the symbol $\sum$ is purely formal. 
    
    \begin{proposition}\label{PropositionHahnProduct}
        Let $\mathcal{H}$ be an RV-sort. Then the Hahn field $(\mathcal{H}^{(\Gamma)},+,\cdot,0,1,\val)$ is a henselian valued field. In addition, we have that
        \[
        \mathcal{H}_{\{0\}}(\mathcal{H}^{(\Gamma)})\simeq\mathcal{H}
        \]
    \end{proposition}
    
    \begin{proof}
        As for the Hahn field of fields with respect of an ordered group, the difficult part is to show that every non-zero element of $\mathcal{H}^{(\Gamma)}$ has a multiplicative inverse. We first show that it is a spherically complete ring. Then we will deduce that it is indeed a field.  Let $a=\sum\limits_{\delta \in \Gamma} \mathbf{a}_\delta$, $  b=\sum\limits_{\epsilon \in \Gamma} \mathbf{b}_\epsilon$,  $c=\sum\limits_{\zeta \in \Gamma}\mathbf{c}_\zeta$ be three elements of $\mathcal{H}^{(\Gamma)}$. 
        \begin{itemize}
        
            \item (associativity for $+$) If $\textbf{a},\textbf{b},\textbf{c} \in \mathcal{H}$ with $\val_{\mathcal{H}}(\textbf{a})=\val_{\mathcal{H}}(\textbf{b})=\val_{\mathcal{H}}(\textbf{c})$, then $(\textbf{a}\oplus \textbf{b}) \oplus \textbf{c} = \textbf{a} \oplus (\textbf{b} \oplus \textbf{c})$. Then associativity for $+$ in  $\mathcal{H}^{(\Gamma)}$ is clear as we sum componentwise.
            \item (commutativity for $+$) Again this is  clear as $\oplus$ is commutative in $\mathcal{H}$.
            \item (neutral element for $+$): $0 :=\sum_\gamma \mathbf{0} \in \mathcal{H}^{(\Gamma)}$ is a neutral element, as $ \mathbf{0} \in \mathcal{H}$ is a neutral element for $\oplus$.
            \item (inverse for $+$) If $a=\sum\limits_{\gamma\in \Gamma}\mathbf{a}_\gamma$, the inverse of $a$ is given by $a=\sum\limits_{\gamma\in \Gamma}-\mathbf{a}_\gamma$. The support being the same, it is an element of $\mathcal{H}^{(\Gamma)}$. 
        \end{itemize}
        
        \begin{itemize}
            \item (associativity for $\cdot$) As the multiplication in $\mathcal{H}$ is associative, a simple computation gives: 
            \[(a\cdot b) \cdot c = \bigoplus\limits_{\delta+\epsilon+\zeta =\gamma} (\mathbf{a}_\delta \cdot \mathbf{b}_\epsilon)\cdot \mathbf{c}_\zeta =\bigoplus\limits_{\delta+\epsilon+\zeta =\gamma} \mathbf{a}_\delta \cdot (\mathbf{b}_\epsilon\cdot \mathbf{c}_\zeta) = a\cdot (b \cdot c).\]
            
            \item (commutativity for $\cdot$) Similarly, this follows from the commutativity of the multiplication in $\mathcal{H}$.
        
            \item (1 is a neutral element for $\cdot$) 
            Immediate as $\textbf{1}\in \mathcal{H}$ is a neutral element for the multiplication in $\mathcal{H}$.
            \item (distributivity)  We have
        \begin{eqnarray*} (a+b)\cdot c &=& \sum_{\gamma \in \Gamma} \bigoplus_{\delta+\epsilon = \gamma} (\mathbf{a}_\delta \oplus \mathbf{b}_\delta) \cdot \mathbf{c}_\epsilon = \sum_{\gamma \in \Gamma} \bigoplus_{\delta+\epsilon = \gamma} (\mathbf{a}_\delta\cdot \mathbf{c}_\epsilon \oplus \mathbf{b}_\delta\cdot \mathbf{c}_\epsilon) \\
        &=& \sum_{\gamma \in \Gamma} \bigoplus_{\delta+\epsilon = \gamma} \mathbf{a}_\delta\cdot \mathbf{c}_\epsilon + \sum_{\gamma \in \Gamma} \bigoplus_{\delta+\epsilon = \gamma} \mathbf{b}_\delta\cdot \mathbf{c}_\epsilon\\
        &=&a\cdot c+b\cdot c.
        \end{eqnarray*}
            \item (valuation) Clearly, $\val$ is a homomorphism of groups. The ultrametric inequality also follows immediately from the definition. 
            \item (spherically complete) We give here a usual diagonal argument. Let $(a^i)_{i<\lambda}$ be a pseudo-Cauchy sequence in $\mathcal{H}^{(\Gamma)}$, where $\lambda$ is any limit ordinal. There is $i_0$ such that for all $i_0<i<j<k$, $\val_{\mathcal{H}}(a^i-a^j)<\val_{\mathcal{H}}(a^j-a^k)$. For $i>i_0$, we denote by $\gamma_i$ the value $\val_{\mathcal{H}}(a^i-a^{i+1})$.  We define 
            \[\mathbf{a}_\gamma :=  \begin{cases} \mathbf{a}_\gamma^i \text{ if }\gamma_i>\gamma \text{ for some  }i<\lambda,\\
        0 \text{ otherwise. } 
        \end{cases}\]
        By definition of the valuation, $\mathbf{a}_\gamma^i$ does not depend on the choice of $i>i_0$ such that $\gamma_i>\gamma$.  Let $a = \sum_\gamma \mathbf{a}_\gamma$, we get for all $i\in I$ that $\val(a - a^i)>\gamma_i$. Then $a$ is a pseudo limit of $(a^i)_i$. We have proved that any pseudo-Cauchy sequence admits a pseudo-limit. 
            \item (multiplicative inverse) Let $a=\sum_\gamma \mathbf{a}_\gamma\in \mathcal{H}^{(\Gamma)}$. Assume it has no inverse and consider the set 
            \[\Delta :=\lbrace \val(1-a\cdot b) \ \vert \ b \in \mathcal{H}^{(\Gamma)} \rbrace.\]
            This set has no maximal element. Indeed, notice first that if $\gamma=\val(1-a\cdot b)$, then $\gamma<\infty$ as $a$ has no inverse. It follows that if $\mathbf{c}_\gamma \in \mathcal{H}$ is the coefficient of value $\gamma$ in $c=1-a\cdot b$, we have $\val(1-a\cdot(b-\mathbf{c}_\gamma\cdot \mathbf{a}_{\val(a)}^{-1}))>\gamma$. So $\Delta$ has no maximal element.
            Let $(\val(1-a\cdot b_\nu))_{\nu\in \lambda}$ be an co-final increasing sequence in $\Delta$. Then, $\lambda$ is a limit ordinal. By definition, $(a\cdot b_\nu)_\nu$ is pseudo-Cauchy with pseudo-limit $1$.  Then, $(b_\nu)_{\nu\in \lambda}$ is a pseudo-Cauchy sequence (as the multiplication by $a$ preserves pseudo-Cauchy sequences). It converges to an element $b\in \mathcal{H}^{(\Gamma)}$, which satisfies $\val(1-a\cdot b) > \Delta$. Contradiction.
        \end{itemize}
        We have proved that $\mathcal{H}^{(\Gamma)}$ is a spherically complete valued field, so in particular it is henselian. Every element $\sum_\gamma \mathbf{a}_\gamma \in \mathcal{H}^{(\Gamma)}$ can be written as $ {\mathbf{a}_\delta(1+\sum_\gamma \mathbf{a}_\gamma/\mathbf{a}_\delta)}$ with $\val(\sum_\gamma \textbf{a}_\gamma/\mathbf{a}_\delta))>0$. Then, one can easily verify that the map $a (1+\mathfrak{m}_{\{0\}}) \mapsto a_\gamma$ gives an isomorphism of hyperfields $\mathcal{H}_{\{0\}}(\mathcal{H}^{(\Gamma)}) \simeq\mathcal{H}$.
    \end{proof}

\section{The main result}\label{SectionMainResult}

The next theorem answers Question \ref{question}. The proofs is built from Proposition  \ref{PropositionInverseLimitHyperfields} and Proposition \ref{PropositionHahnProduct}.
\begin{theorem}
    Let $\Gamma$ be an ordered abelian group. Assume that $\Delta$ is a convex subgroup of $\Gamma$ and denote by $\mathcal{N}:= \Delta_{\geq 0}$ its positive part.
    Let $(\mathcal{H}_i)_{i\in I}$ be an isometric system of valued hyperfields $\mathcal{H}_i$ with value group $\Gamma$ and norm $\mathcal{N}_i\subseteq\mathcal{N}$ such that $\mathcal{N}=\cup_i \mathcal{N}_i$. Then there is a henselian valued field $\mathcal{K}$ with value group $\Gamma$ and such that $\mathcal{H}_{\mathcal{N}_i}(\mathcal{K})\simeq\mathcal{H}_i$ for every $i\in I$.

\end{theorem}

\begin{proof}
    To ease the notation and thanks to Proposition \ref{PropositionQuotientHyperfields} and Lemma \ref{LemmaIsometryIndicesIsomophism}, we can complete the sequence  
     $(\mathcal{H}_i)_{i\in I}$ into a sequence $(\mathcal{H}_\rho)_{\rho\subset \mathcal{N}}$ where $\mathcal{N}=\cup_{i\in I} \mathcal{N}_i$. Consider  the inverse limit $\mathcal{H}$. By Proposition \ref{PropositionInverseLimitHyperfields},  $\mathcal{H}$ is a stringent valued hyperfield with a valuation $\val_{\mathcal{H}}: \mathcal{H} \rightarrow \Gamma$. Consider the equivalence relation $\sim$ on $\mathcal{H}$ as defined in Section \ref{SectionAxiomatisationStringentHyperfields}. Since
     \[ 
     a\sim b \ \overset{def}{\Longleftrightarrow} \  a+b \neq \{a\} \wedge a+b \neq \{b\} \ \Longleftrightarrow \  \val_{\mathcal{H}}(a)-\val_{\mathcal{H}}(b)\in \Delta,
     \]
     we have that $\mathcal{H}^\times/{\sim} \simeq \Gamma/\Delta$.
     It follows that, as a stringent hyperfield, $\mathcal{H}$ arises from the following sequence-structure:
        \[
        \begin{tikzcd}
        \{1\}\arrow[r]& \mathcal{F}^\times \arrow[r,"\iota"]& \mathcal{H}^\times\arrow[r,"w'"]& \Gamma/\Delta\arrow[r]& \{0\}.
        \end{tikzcd}
        \]
    where $\mathcal{F}$ is a certain field and $w':\mathcal{H} \rightarrow \Gamma/\Delta$ is the valuation on $\mathcal{H}$ defined as $w'(\mathbf{a}):=\val_\mathcal{H}(\rv(a))+\Delta$.
    
    Using the Hahn construction of Proposition \ref{PropositionHahnProduct}, we obtain now a henselian valued field
    $\mathcal{K}$ with valuation $\val':\mathcal{K} \rightarrow \Gamma/\Delta$ and projection $\rv: \mathcal{K} \rightarrow \mathcal{H}$ such that $w'(\rv(a))=\val'(a)$ for all $a\in \mathcal{K}$.  
     For $a\in \mathcal{K}$, we define $\val(a):= \val_{\mathcal{H}}(\rv(a))$. We now show that $\val$ is a valuation on $\mathcal{K}$:
    \begin{itemize}
        \item $\val(a)= \val_{\mathcal{H}}(\rv(a))=\infty$ if and only if $\rv(a)= 0$ if and only if $a=0$, for all $a\in\mathcal{K}$.
        \item $\val(a b)= \val_{\mathcal{H}}(\rv(ab))= \val_{\mathcal{H}}(\rv(a)\rv(b))=\val_{\mathcal{H}}(\rv(a))+\val_{\mathcal{H}}(\rv(b))$, for all $a,b\in \mathcal{K}$.
        \item For $a,b\in \mathcal{K}$, let us show that $\val(a+b)\geq \min\{\val(a),\val(b)\}$. If $\rv(a) \neq -\rv(b)$, then $\rv(a+b)$ is the unique element in $\rv(a)+\rv(b)$ (since $\mathcal{H}$ is stringent and  $\rv(a+b) \in \rv(a)+\rv(b)$ by Remark \ref{RemarkClassofSumIsAlwaysInSumClass}) and
    $\val(a+b)= \val_{\mathcal{H}}(\rv(a) + \rv(b)) \geq \min(\val_{\mathcal{H}}(\rv(a)), \val_{\mathcal{H}}(\rv(b)))$
    since $\val_\mathcal{H}$ is a valuation on $\mathcal{H}$. On the other hand, if $\rv(a) = -\rv(b)$, then in particular $w'(\rv(a+b))>w'(\rv(a)
    )$ by the norm axiom applied to $w'$. Since $w'(rv(a))=\val(a)+\Delta$, we deduce that 
    \[
    \val(a+b)= \val_{\mathcal{H}}(\rv(a+b)) > \val_\mathcal{H}(\rv(a))=\min(\val_{\mathcal{H}}(\rv(a)), \val_{\mathcal{H}}(\rv(b))). 
    \]
    \end{itemize}

We can now prove that $\mathcal{H}_\rho(\mathcal{K})\simeq\mathcal{H}_\rho$ for every $\rho\subseteq\mathcal{N}$. Since 
\[
\mathfrak{m}_\mathcal{N}:= \{x\in \mathcal{K} \ \vert \ \val(x)>\mathcal{N}  \} = \{x\in \mathcal{K} \ \vert \ w'(x)>0\},
\] 
we have by Proposition \ref{PropositionHahnProduct} that $\mathcal{H}_{\mathcal{N}}(\mathcal{K})= \mathcal{K}_{1+\mathfrak{m}_\mathcal{N}}$ is isomorphic to $\mathcal{H}$. Moreover, by construction, the valuation induced by $\val$ on $\mathcal{H}$ is exactly $\val_\mathcal{H}$. As $\mathcal{H}$ is the inverse limit of $(\mathcal{H}_\rho)_{\rho \subseteq \mathcal{N}}$, it follows from Proposition \  \ref{PropositionInverseLimitHyperfields} \ (3) that
    \[
    \mathcal{H}_{\rho}(\mathcal{K})\simeq \mathcal{H}_{\mathcal{N}}(\mathcal{K})_{1+\mathfrak{m}_\rho} \simeq \mathcal{H}_{\rho}
    \]
    for all initial segments $\rho \subseteq \mathcal{N}$.  This completes the proof.
\end{proof}

\begin{example}
We consider the field $\mathcal{K}=\mathbb{Q}(x)(y)$ equipped with the valuation $\val:= \val_x \circ \val_y : \mathcal{K} \rightarrow \mathbb{Z}\times \mathbb{Z}$ of the $x$-adic and $y$-adic valuation ($\val(x)=(0,1)<\val(y)=(1,0)$) as in Example \ref{ExamplesProjectiveLimiteHyperfields}.\par
If $\mathcal{H}:=\varprojlim_{n}\mathcal{H}_{(0,n)}(\mathcal{K})$, then the Hahn product $\mathcal{H}^{(\mathbb{Z})}$ is isomorphic to $\mathbb{Q}((x))((y))$. 

 Similarly, if $\mathcal{K}'=\mathbb{Q}(x)(y)(z)$, and $\mathcal{H}'$ is the inverse limit $\varprojlim_{n}\mathcal{H}_{(0,0,n)}(\mathcal{K'})$ then the Hahn product $\mathcal{H'}^{(\mathbb{Z}\times \mathbb{Z})}$ gives  $\mathbb{Q}((x))(y)((z))$, i.e. the completion in $x$ and $z$  but not in $y$. If $\mathcal{H}''$ is the inverse limit $\varprojlim_{n}\mathcal{H}_{(0,n,0)}(\mathcal{K''})=\varprojlim_{n}\mathcal{H}_{(0,n,n)}(\mathcal{K''})$ then the Hahn product $\mathcal{H''}^{(\mathbb{Z})}$ gives $\mathbb{Q}(x)((y))((z))$, i.e. the completion in $y$ and $z$ but not in $x$.
\end{example}

We conclude with a question which, to our knowledge, remains open:

\begin{question}
    Are all valued hyperfields with value group $\Gamma$ of the form $\mathcal{H}_\rho(\mathcal{K})$ for some valued field $\mathcal{K}$ and some initial segment $\rho$ of $\Gamma_    {\geq 0}$?
\end{question}

    \bibliographystyle{plain}
    \bibliography{main}
    
\end{document}